\documentclass[11pt]{article}
\usepackage{latexsym,amsmath,amssymb,graphics,amscd}
\makeatletter
\@addtoreset{figure}{section}
\def\thefigure{\thesection.\@arabic\c@figure}
\def\fps@figure{h,t}
\@addtoreset{table}{bsection}

\def\thetable{\thesection.\@arabic\c@table}
\def\fps@table{h, t}
\@addtoreset{equation}{section}

\makeatother

\begin{document}

\newtheorem{theorem}{Theorem}[section]
\newtheorem{definition}[theorem]{Definition}
\newtheorem{lemma}[theorem]{Lemma}
\newtheorem{remark}[theorem]{Remark}
\newtheorem{proposition}[theorem]{Proposition}
\newtheorem{corollary}[theorem]{Corollary}
\newtheorem{example}[theorem]{Example}

\newcommand{\bfi}{\bfseries\itshape}
\newcommand{\mmu}{\ensuremath{\mathbf{J}^{-1}(\mu)/G_{\mu}}}
\newcommand{\mmuh}{\ensuremath{(\mathbf{J}^{-1}(\mu)\cap
M_{(H)}^{G_{\mu}})/G_{\mu}}}
\newcommand{\nmmuh}{\ensuremath{\mathbf{J}^{-1}(\mu)\cap M_{(H)}^{G_{\mu}}}}
\newcommand{\momu}{\ensuremath{\mathbf{J}^{-1}(\mathcal{O}_{\mu})/G}}
\newcommand{\momuh}{\ensuremath{\mathbf{J}^{-1}(\mathcal{O}_{\mu})\cap 
M_{(H)}/G}}
\newcommand{\nmomuh}{\ensuremath{\mathbf{J}^{-1}(\mathcal{O}_{\mu})\cap 
M_{(H)}}}
\newcommand{\omu}{\ensuremath{\mathcal{O}_{\mu}}}
\newcommand{\ozero}{\ensuremath{\widehat{\mathcal{O}}_{0}}}
\newcommand{\oalpha}{\ensuremath{\mathcal{O}}_{\alpha}}
\newcommand{\ddto}{\ensuremath{\left.\frac{d}{dt}\right|_{t=0}}}
\newcommand{\ddt}{\ensuremath{\frac{d}{dt}}}
\newcommand{\ddso}{\ensuremath{\left.\frac{d}{ds}\right|_{s=0}}}
\newcommand{\sm}{\ensuremath{_{\mu}}}
\newcommand{\seta}{\ensuremath{_{\eta}}}
\newcommand{\sbeta}{\ensuremath{_{\beta}}}
\newcommand{\sele}{\ensuremath{_{L}}}
\newcommand{\seme}{\ensuremath{_{m}}}
\newcommand{\ele}{\ensuremath{N(H)/H}}
\newcommand{\ene}{\ensuremath{\mathfrak{n}}}
\newcommand{\enesmu}{\ensuremath{\mathfrak{n}_{\mu}}}
\newcommand{\gele}{\ensuremath{\mathfrak{l}}}
\newcommand{\elesmu}{\ensuremath{\mathfrak{l}_{\mu}}}
\newcommand{\snu}{\ensuremath{_{\nu}}}
\newcommand{\sem}{\ensuremath{_{m}}}
\newcommand{\smu}{\ensuremath{_{\mu}^{(H)}}}
\newcommand{\smus}{\ensuremath{_{\mu}^{(H) *}}}
\newcommand{\som}{\ensuremath{_{\mathcal{O}_{\mu}}}}
\newcommand{\somu}{\ensuremath{_{\mathcal{O}_{\mu}}^{(H)}}}
\newcommand{\somus}{\ensuremath{_{\mathcal{O}_{\mu}}^{(H) *}}}
\newcommand{\J}{\ensuremath{\mathbf{J}}}
\newcommand{\K}{\ensuremath{\mathbf{K}}}
\newcommand{\g}{\ensuremath{\mathfrak{g}}}
\newcommand{\gt}{\ensuremath{\mathfrak{t}}}
\newcommand{\gts}{\ensuremath{\mathfrak{t}^{*}}}
\newcommand{\sgt}{\ensuremath{_{\mathfrak{t}}}}
\newcommand{\gmu}{\ensuremath{\mathfrak{g}_{\mu}}}
\newcommand{\h}{\ensuremath{\mathfrak{h}}}
\newcommand{\m}{\ensuremath{\mathfrak{m}}}
\newcommand{\s}{\ensuremath{\mathfrak{s}}}
\newcommand{\lamba}{\ensuremath{\bar{\lambda}}}
\newcommand{\eba}{\ensuremath{\bar{\eta}}}
\newcommand{\hba}{\ensuremath{\bar{h}}}
\newcommand{\xiba}{\ensuremath{\bar{\xi}}}
\newcommand{\rhoba}{\ensuremath{\bar{\rho}}}
\newcommand{\tauba}{\ensuremath{\bar{\tau}}}
\newcommand{\pe}{\ensuremath{\mathfrak{p}}}
\newcommand{\q}{\ensuremath{\mathfrak{q}}}
\newcommand{\qs}{\ensuremath{\mathfrak{q}^{*}}}
\newcommand{\subh}{\ensuremath{_{\mathfrak{h}}}}
\newcommand{\subm}{\ensuremath{_{\mathfrak{m}}}}
\newcommand{\subs}{\ensuremath{_{\mathfrak{s}}}}
\newcommand{\gs}{\ensuremath{\mathfrak{g}^{*}}}
\newcommand{\gmus}{\ensuremath{\mathfrak{g}_{\mu}^{*}}}
\newcommand{\hs}{\ensuremath{\mathfrak{h}^{*}}}
\newcommand{\ms}{\ensuremath{\mathfrak{m}^{*}}}
\newcommand{\subms}{\ensuremath{_{\mathfrak{m}^{*}}}}
\newcommand{\subv}{\ensuremath{_{V}}}
\newcommand{\bd}{\ensuremath{\mathbf{d}}}
\newcommand{\subt}{\ensuremath{_{t}}}
\newcommand{\sss}{\ensuremath{\mathfrak{s}^{*}}}
\newcommand{\inv}{\ensuremath{^{-1}}}
\newcommand{\suh}{\ensuremath{^{H}}}
\newcommand{\ssuh}{\ensuremath{^{(H)}}}
\newcommand{\sus}{\ensuremath{^{*}}}
\newcommand{\suxi}{\ensuremath{^{\xi}}}
\newcommand{\sucero}{\ensuremath{^{\circ}}}
\newcommand{\scero}{\ensuremath{_{\circ}}}
\newcommand{\szero}{\ensuremath{_{0}}}
\newcommand{\slo}{\ensuremath{_{\lambda_{\circ}}}}
\newcommand{\lo}{\ensuremath{\lambda_{\circ}}}
\newcommand{\no}{\ensuremath{\nu_{\circ}}}
\newcommand{\sno}{\ensuremath{_{\nu_{\circ}}}}
\newcommand{\sh}{\ensuremath{_{H}}}
\newcommand{\ssh}{\ensuremath{_{(H)}}}
\newcommand{\suinf}{\ensuremath{^{\infty}}}
\newcommand{\ngmu}{\ensuremath{N_{G_{\mu}}(H)}}
\newcommand{\ngmuh}{\ensuremath{N_{G_{\mu}}(H)/H}}
\newcommand{\ad}{\ensuremath{{\rm ad}}}
\newcommand{\Ad}{\ensuremath{{\rm Ad}}}
\newcommand{\soh}{\ensuremath{_{0}^{(H)}}}
\newcommand{\sok}{\ensuremath{_{0}^{(K)}}}

\makeatletter
\title{Hamiltonian Hopf bifurcation with symmetry}
\author{Pascal Chossat$^{1}$, Juan-Pablo Ortega$^{2}$, and Tudor S. Ratiu$^{3}$}
\addtocounter{footnote}{1}
\footnotetext{Institut Nonlin\'eaire de Nice, UMR 129, CNRS-UNSA, 1361, route des Lucioles, 06560 Valbonne, France. \texttt{chossat@inln.cnrs.fr.}}
\addtocounter{footnote}{1}
\footnotetext{D\'epartement de Math\'ematiques, 
\'Ecole Polytechnique F\'ed\'erale de Lausanne. CH--1015 Lausanne. Switzerland.
\texttt{Juan-Pablo.Ortega@epfl.ch}.}
\addtocounter{footnote}{1}
\footnotetext{Department of Mathematics, 
University of California, Santa Cruz, Santa Cruz, CA 95064, USA, and D\'epartement de Math\'ematiques, 
\'Ecole Polytechnique F\'ed\'erale de Lausanne. CH--1015 Lausanne. Switzerland.
\texttt{Tudor.Ratiu@epfl.ch}. Research partially supported by
NSF Grant DMS-9802378  and FNS Grant 21-54138.98.}
\date{November 27, 2000}
\maketitle
\makeatother

\begin{abstract}
In this paper we study the appearance of branches of relative periodic orbits in Hamiltonian Hopf bifurcation processes in the presence of compact symmetry groups that do not generically exist in the dissipative framework. The theoretical study is illustrated with several examples.
\end{abstract}

\section{Introduction}

\begin{figure}[htb]
\begin{center}
\includegraphics{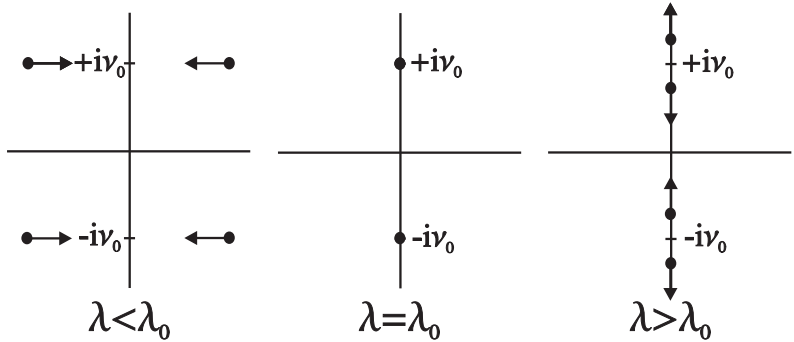}
\caption{Motion of eigenvalues in a Hamiltonian Hopf bifurcation.}\label{fig:hopf}
\end{center}
\end{figure}

Let $(V, \omega)$ be a symplectic vector space and $G$ be a compact Lie group acting linearly and symplectically on $V$. Let now be a one--parameter family of $G$--invariant Hamiltonians $h_\lambda\in C\suinf(V)^G$ such that for each value of the parameter $\lambda$, the origin is an equilibrium of the associated Hamiltonian vector field, that is, $\bd h_\lambda(0)=0$ for arbitrary $\lambda$. In this paper we will study the nonlinear implications of the following linear behavior: suppose that there is a value of the parameter $\lo$ and a pair of  eigenvalues $\pm i\no$ in the spectrum of the  linearization at zero of the dynamics induced by the Hamiltonian vector field $X_{h\slo}$ that behave as in Figure~\ref{fig:hopf} when we move the parameter $\lambda$ around $\lo$.  Such a behavior in the parametrical motion of the eigenvalues is usually referred to as {\bfi Hamiltonian Hopf bifurcation}~\cite{van der meer 85}, denomination that we will use here, even though it also appears in the literature as $1:-1$ resonance, $1:1$ non--semisimple resonance, and Krein collision. The reference to the Hopf bifurcation comes from the analogy with the codimension one non--conservative case in which a one-parameter family of vector fields has a pair of eigenvalues that cross the imaginary axis at a critical value of the parameter (the "classical" Hopf bifurcation). The case of $G$-equivariant vector fields ($G$ compact) has led to the successful theory of Hopf bifurcation with symmetry which was initiated by \cite{g1} and which was exposed in its most achieved form in \cite{field} (see also \cite{cl} for a comprehensive exposition).

The history of the Hamiltonian Hopf bifurcation in the non-symmetric setup is very long and we shall not attempt to survey it here. We just refer to~\cite{meyer schmidt, meyer 86, van der meer 85, van der meer 96, bridges 90, constrained paper} and references therein for discussions.

The only works that we know of dealing with the Hamiltonian symmetric case are~\cite{van der meer 90} and~\cite{rhombic 96}. In the first paper  it is shown that one can apply the non--symmetric results  on Hamiltonian Hopf bifurcation on some of the fixed point spaces corresponding to isotropy subgroups of the symmetry of the system, provided that certain dimensional restrictions are fulfilled. \cite{rhombic 96} studies branches of (stable) three--tori that can be obtained out of a Hamiltonian Hopf bifurcation process with a symmetry given by the semi--direct product of $D_2$ with $T^2\times S^1$. See also~\cite{bridges 90a}.

Natural dynamical elements that show up in the study of systems that present a continuous symmetry group $G$ are the so called {\bfi relative equilibria (RE)} and {\bfi relative periodic orbits (RPOs)}, that is,  motions that project onto equilibria and periodic orbits in the quotient space $V/G$, respectively. In our work we will see that whenever a Hopf--like motion of eigenvalues occurs in a Hamiltonian system with symmetry, one can prove the existence of periodic and relative periodic motions at the non linear level, for values of the parameter nearby $\lo$, whose number we will estimate at each energy level.  The existence of periodic motions, in the presence of some dimensional restrictions that we have eliminated, was something known to the above quoted authors. As to the relative periodic orbits, they are found in those papers only after a reduction has been performed that makes the problem equivalent to that of searching periodic orbits in the relevant quotient space. Since this reduction cannot always be carried out in a straightforward manner we will follow an approach in which the existence of relative periodic orbits is proved in the original space $V$. 

\medskip

Our approach to this problem will be based in the combined use of five tools:
\begin{description}
\item[(i)] Reduction method of Vanderbauwhede and van der Meer~\cite{vanderbauwhede 95}: it allows us to substitute the search of periodic and relative periodic orbits by the search of relative equilibria of a $S^1$--symmetric associated Hamiltonian system (usually referred to as {\bfi normal form}).
\item[(ii)] Generic structure of the generalized eigenspaces corresponding to the colliding eigenvalues~\cite{dellnitz melbourne marsden}: it determines the most plausible reduced space in which we should work after applying~\textbf{(i)}.
\item[(iii)] Equivariant Williamson normal form~\cite{melbourne dellnitz 93}: it is used to normalize  the linear term of the equation that defines the $S^1$--relative equilibria that we are looking for.
\item[(iv)] Lyapunov-Schmidt reduction of the finite dimensional equation that defines the $S^1$--relative equilibria and formulation of the problem in terms of a  bifurcation equation of gradient nature~\cite{pascal} with very specific equivariance properties.
\item[(v)] Solution of the bifurcation equation using either topological or analytical methods.
\end{description} 

The paper is structured as follows: 
\begin{itemize}
\item In Section~\ref{preliminaries} we briefly review the abovementioned tools, set the notation that will be used throughout the paper, and explain in detail the hypotheses under which we will work, along with their implications. The expert can skip this section and use it just as a glossary. 
\item Section~\ref{the main result} is devoted to  Theorem~\ref{relative periodic theorem}, which provides a lower estimate on the number of periodic and relative periodic branches that bifurcate from the origin if there is a collision of eigenvalues as in Figure~\ref{fig:hopf}. 
\item In Section~\ref{example} we study a system of two nonlinearly coupled harmonic oscillators in the presence of a magnetic field, that will lead us to the consideration of the general case of the Hamiltonian Hopf bifurcation in the presence of a $O(2)$ symmetry. We will see that, in contrast with the dissipative case, the $O(2)$--symmetry in a Hamiltonian Hopf bifurcation process gives rise to the appearance of numerous relative periodic motions. This example will also show that, in general, the topological methods utilized in Theorem~\ref{relative periodic theorem} are not powerful enough to detect all the periodic and relative periodic elements of a particular system with a given symmetry, that is, the generality of this result is paid with its lack  of sharpness. 
This circumstance will motivate a more hands--on approach to the problem in Section~\ref{manual treatment}, where we will see that under additional hypotheses on the group action, sharper general results can be formulated that give account of all the dynamical richness evidenced in the example in Section~\ref{example}. 
\item In Section~\ref{spheric example} we use the previous results to show the existence of RPOs in Hamiltonian Hopf phenomena with spheric symmetry.
\item For the sake of the clarity in the exposition, the proofs of some of the technical results needed in the main theorem are relegated to an appendix (Section~\ref{appendix}) at the end of the paper.
\end{itemize}

\medskip

\section{Preliminaries and setup}
\label{preliminaries}

All throughout this paper, our discussions will mostly take place in a finite dimensional symplectic vector space $(V, \omega)$ on which the compact Lie group $G$ acts linearly and canonically, that is, respecting the symplectic form $\omega$. We will be interested in a one--parameter family of $G$--equivariant Hamiltonian vector fields $X_{h_\lambda}$, induced by the  family of $G$--invariant Hamiltonians $h_\lambda\in C\suinf(V)^G$,  $\lambda\in\mathbb{R}$, such that:
\begin{description}
\item[(H1)] $h_\lambda(0)=0$ and $\bd h_\lambda(0)=0$ for all $\lambda$.
\item[(H2)] There is a value $\lo$ of the parameter $\lambda$ for which the $G$--equivariant infinitesimally symplectic linear map $A\slo:=D\subv X_{h\slo}(0)$ is non singular and has $\pm i\no$ in its spectrum ($\no\neq 0$).
\end{description}

\medskip

In the following paragraphs we introduce some tools and notations that will be used all throughout the paper.

\medskip

\noindent\textbf{The resonance space} Let 
$(V, \omega)$ be a symplectic vector space. 
It is easy to show that there is a bijection between 
linear Hamiltonian vector fields on $(V, \omega)$ 
and quadratic forms on $V$. Indeed, if 
$A:V\rightarrow V$ is an infinitesimally 
symplectic linear map, that is, a linear Hamiltonian 
vector field on $(V, \omega)$, its corresponding 
Hamiltonian function is given by
\[
Q_A(v):=\frac{1}{2}\omega(Av, v),
\qquad\text{for any $v\in V$.}
\]
Also, since $A$ belongs to the symplectic Lie algebra $\mathfrak{sp}(V)$, it admits a unique {\bfi Jordan--Chevalley decomposition}~\cite{humphreys, vanderbauwhede 95} of the form $A=A_s+A_n$, where $A_s\in \mathfrak{sp}(V)$ is semisimple (complex diagonalizable), $A_n\in \mathfrak{sp}(V)$ is nilpotent, and $[A_s, A_n]=0$.
If
$i\no$ is one of the eigenvalues of $A\in \mathfrak{sp}(V)$ and 
$T\sno:=\frac{2\pi}{\no}$, we define the 
{\bfi resonance space} $U\sno$ of $A$ with 
{\bfi primitive period} $T\sno$ as 
\[
U\sno:=\ker\left(e^{A_sT\sno}-I\right).
\]
The resonance space $U\sno$ has the 
following properties (see \cite{williamson, 
golubitsky stewart 87, vanderbauwhede 95}):
\begin{description}
\item[(i)] $U\sno$ is equal to the direct sum 
of the real generalized eigenspaces of $A$ 
corresponding to eigenvalues of the form 
$\pm ik\no$, with $k\in\mathbb{N}\sus$.
\item[(ii)] The pair $(U\sno, \omega|_{U\sno})$ 
is a symplectic subspace of $(V, \omega)$.
\item[(iii)] The mapping $\theta\in S^1\mapsto 
e^{\frac{\theta}{\no} A_s}|_{U\sno}$ generates a symplectic 
$S^1$--linear action  on $(U\sno, \omega|_{U\sno})$, whose associated equivariant momentum map will be denoted by  $\J:U\sno\rightarrow{\rm Lie}(S^1)\sus\simeq\mathbb{R}$. 
\item[(iv)] If $(V, \omega)$ is a symplectic 
representation space of the Lie group $G$ 
and the Hamiltonian vector field $A$ is 
$G$--equivariant (equivalently, the quadratic 
form $Q_A$ is $G$--invariant), then the 
symplectic resonance subspace 
$(U\sno, \omega|_{U\sno})$ is also 
$G$--invariant (this follows from the uniqueness of the Jordan--Chevalley decomposition of $A$, which implies that if $A$ is $G$--equivariant, so is $A_s$). Moreover, the $S^1$ 
and $G$ actions on $(U\sno, \omega|_{U\sno})$ 
commute, which therefore defines a  symplectic 
linear action of $G\times S^1$ on $U_{\no}$. See the Appendix (Section~\ref{appendix}) for a sketch of the proof of some of these facts.
\end{description}

\medskip

\noindent\textbf{The normal form reduction} \cite{van der meer 85, van der meer 90, vanderbauwhede 95} Let $(V, \omega, h_\lambda)$ be a $\lambda$--parameter family ($\lambda\in \Lambda$, where $\Lambda$ is a Banach space) of $G$--Hamiltonian systems such that $h\slo (0)=0$, $\bd h\slo (0)=0$, and the $G$--equivariant infinitesimally symplectic linear map $A:=DX_{h\slo}(0)$ is non singular and has $\pm i\no$ as eigenvalues. Let $(U\sno, \omega|_{U\sno})$ be the resonance space of $A$ with primitive period $T\sno$. For each $k\geq 0$ there are a $C^k$--mapping $\psi:U\sno\times\Lambda\rightarrow V$ and a $C^{k+1}$--mapping $\widehat{h_\lambda}:U\sno\times\Lambda\rightarrow\mathbb{R}$ such that  $\psi(0,\lambda)=0$, for all $\lambda\in\Lambda$,   $D_{U\sno}\psi(0,\lo)=\mathbb{I}_{U\sno}$,  and $\widehat{h_\lambda}$ is a $G\times S^1$--invariant function  that coincides with $h_\lambda$ up to order $k+1$. 
The interest of normalization is given by the fact that one can prove~\cite[Theorem 3.2]{vanderbauwhede 95} that if we stay close enough to zero in $U\sno$ and to $\lo\in \Lambda$, then the  $S^1$--relative equilibria  of the $G\times S^1$--invariant Hamiltonian $\widehat{h_\lambda}$ are mapped by $\psi(\cdot,\lambda)$ to the set of periodic solutions of $(V, \omega, h_\lambda)$ in a neighborhood of $0\in V$ with periods close to $T\sno$. Hence, in our future discussion we will substitute the problem of searching periodic orbits for $(V, \omega, h_\lambda)$ by that of searching the $S^1$--relative equilibria of the $G\times S^1$--invariant family of Hamiltonian systems $(U\sno, \omega|_{U\sno}, \widehat{h_\lambda})$, that will be referred to as the {\bfi equivalent system}. Note that the properties of  $\psi$ imply that
\begin{equation}
\label{restriction to resonance}
\mathcal{A}:=A|_{U\sno}=D\subv X_{h\slo}(0)|_{U\sno}=D\subv X_{h\slo|_{U\sno}}(0)=D\subv X_{\widehat{h\slo}}(0).
\end{equation}

\medskip

\noindent\textbf{Generic structure of the resonance space $U\sno$ and canonical form of the symplectic pair $(\omega|_{U\sno}, \mathcal{A})$}  After the remarks previously made, we  know that the resonance space $U\sno$ is a $G\times S^1$--symplectic vector space. The decomposition of $U\sno$ into $G\times S^1$--irreducible subspaces that can be generically expected when the eigenvalues behave parametrically as in Figure~\ref{fig:hopf} has been studied in~\cite{dellnitz melbourne marsden}, where the authors concluded (Proposition 6.1 (3)) that the only possibility is 
\begin{equation}
\label{decomposition of resonance space}
U\sno=U_1\oplus U_2,
\end{equation}
where $U_1$ and $U_2$ are complex dual irreducible subspaces of $U\sno$ in the sense of~\cite[Theorem 2.1]{mrs}. In all that follows we will assume that we are in this generic situation.

Once we know the decomposition~(\ref{decomposition of resonance space}) of $U\sno$ into irreducibles, the equivariant version of the Williamson normal form of Melbourne and Dellnitz guarantees~\cite[Table 2]{melbourne dellnitz 93} that there is a basis of the vector space $U\sno$ in which the simultaneous matricial expressions of $\omega|_{U\sno}$ and $\mathcal{A}$ are either
\begin{description}
\item[(i)]\noindent\begin{equation}
\label{canonical form}
\mathcal{A}=
\left(
\begin{array}{cc}
\no \mathbb{J}_{2 n}&\mathbb{I}_{2n}\\
\mathbf{0}&\no \mathbb{J}_{2 n}
\end{array}
\right)\quad\text{and}\quad \omega|_{U\sno}=\mathbb{J}_{4n},\quad\text{or}\quad
\end{equation}
\item[(ii)]\noindent
\begin{equation}
\label{canonical form 2}
\mathcal{A}=
\left(
\begin{array}{cc}
\no \mathbb{J}_{2 n}&\mathbb{I}_{2n}\\
\mathbf{0}&\no \mathbb{J}_{2 n}
\end{array}
\right)\quad\text{and}\quad \omega|_{U\sno}=-\mathbb{J}_{4n},
\end{equation}
\end{description}
where $2n=\dim U_1=\dim U_2$,  $\mathbb{I}_{2 n}$ is the $2n$--dimensional identity matrix, and $\mathbb{J}_{2 n}$ is defined as
\[
\mathbb{J}_{2 n}=
\left(
\begin{array}{cc}
\mathbf{0}&-\mathbb{I}_{n}\\
\mathbb{I}_{n}&\mathbf{0}
\end{array}
\right).
\]
Given that the treatment of cases~\textbf{(i)} and~\textbf{(ii)} is completely analogous, we will focus in all that follows in expression~(\ref{canonical form}). Moreover, whenever our family of $G$--Hamiltonian systems falls into the generic situation described in this paragraph, we will say that it satisfies the condition~\textbf{(H3)}. For clarity and future reference we state this condition explicitely:

\begin{description}
\item[(H3)] The resonance space $U\sno$ corresponding to the eigenvalues $\pm i\no$ splits into two complex dual $G\times S^1$--irreducible subspaces. This condition is generic.
\end{description}
 
\medskip

\noindent\textbf{Relative equilibria and critical points} Being consistent with the notation previously introduced, let $\mathcal{A}=\mathcal{A}_s+\mathcal{A}_n$ be the Jordan--Chevalley decomposition of $\mathcal{A}\in\mathfrak{sp}_G(U\sno)$. We will denote by $\J:U\sno\rightarrow{\rm Lie}(S^1)\sus\simeq\mathbb{R}$ the equivariant momentum map associated to the symplectic 
$S^1$--linear action defined by $(\theta, v)\mapsto 
e^{\frac{\theta}{\no} \mathcal{A}_s}v$, $\theta\in S^1$, $v\in U\sno$. Also, for any $\xi\in {\rm Lie}(S^1)\simeq\mathbb{R}$ and any $v\in U\sno$, we will write $\J\suxi(v):=\J(v)\xi$. The linearity of the action implies that, for any $\xi\in {\rm Lie}(S^1)\simeq\mathbb{R}$ and any $v\in U\sno$, the momentum map $\J$ is uniquely determined by the expression
\[
\J\suxi(v)=\frac{1}{2}\omega|_{U\sno}(\xi\cdot v, v),
\]
where the dot in $\xi\cdot v$ means the associated representation of the Lie algebra ${\rm Lie}(S^1)$ on $U\sno$ through the $S^1$--action. More specifically,
\[
\J(v)=\frac{1}{2\no}\omega|_{U\sno}(\mathcal{A}_s v, v).
\]
For future reference we note that this relation implies that
\begin{equation}
\label{hessian of j}
\bd^2\J(0)(v, w)=\omega|_{U\sno}(\mathcal{A}_s v, w),\qquad\text{for any}\qquad v, w\in U\sno,
\end{equation}
which in the basis used to write~(\ref{canonical form}) admits the following matricial expression:
\begin{equation}
\label{hessian of j matrix}
\bd^2\J(0)=
\left(
\begin{array}{cc}
\mathbf{0}&\mathbb{J}_{2n}\\
-\mathbb{J}_{2n}&\mathbf{0}
\end{array}
\right).
\end{equation}

A very interesting feature of the Hamiltonian framework is that the search for relative equilibria reduces to the determination of the critical points of the so called {\bfi augmented Hamiltonian}~\cite{fom}. In the particular case that we are dealing with, this remark translates into saying that the equivalent system $(U\sno, \omega|_{U\sno}, \widehat{h_\lambda})$ has a $S^1$--relative equilibrium at $v\in U\sno$ (which represents a periodic orbit of the original system $(V, \omega, h_\lambda)$ with period near $T\sno$) if and only if there is an element $\xi\in{\rm Lie}(S^1)$ for which 
\begin{equation}
\label{relative equilibrium condition}
\bd(\widehat{h_\lambda}-\J\suxi)(v)=0.
\end{equation}
Whenever we find a pair $(v, \xi)$ that satisfies~(\ref{relative equilibrium condition}), we will say that $v$ is a relative equilibrium with {\bfi velocity} $\xi$.

Expression~(\ref{relative equilibrium condition}) can be written as a gradient equation, which will be exploited profusely in our subsequent discussion. Indeed, let $\langle\cdot,\cdot\rangle$ be a $G\times S^1$--invariant inner product on $U\sno$ (always available by the compactness of $G\times S^1$). For any $v\in U\sno$, we define the gradient $\nabla_{U\sno}(\widehat{h_\lambda}-\J\suxi)(v)$ as the unique element in $U\sno$, such that for $w\in U\sno$ arbitrary
\[
\bd(\widehat{h_\lambda}-\J\suxi)(v)\cdot w=\langle\nabla_{U\sno}(\widehat{h_\lambda}-\J\suxi)(v), w\rangle.
\]

Also for future reference, we recall that the linearization $A_\lambda=D\subv X_{h_\lambda}(0)$ of $X_{h_\lambda}$ at $0\in V$, is a linear $G$--equivariant Hamiltonian vector field with associated quadratic Hamiltonian function $Q_\lambda$ given by 
\[
Q_\lambda(v)=\frac{1}{2}\bd^2 h_\lambda(0)(v, v),
\]
that is:
\[
{\rm{\mathbf i}}_{A_\lambda}\omega=\bd Q_\lambda.
\]
The restriction $\mathcal{A}$ of $A\slo$ to $U\sno$ is of course also Hamiltonian but in this case, by~(\ref{restriction to resonance}), the associated quadratic Hamiltonian function can be expressed in terms of the Hessian at $0$ of the equivalent Hamiltonian $\widehat{h\slo}$ associated to $h\slo$. Indeed,
\begin{equation}
\label{hamilton for associated}
{\rm{\mathbf i}}_{\mathcal{A}}\omega|_{U\sno}=\bd \mathcal{Q}\slo,
\end{equation}
where, for any $v\in U\sno$,
\begin{equation}
\label{fancy q}
\mathcal{Q}\slo(v)=\frac{1}{2}\bd^2h\slo(0)(v, v)=\frac{1}{2}\bd^2\widehat{h\slo}(0)(v, v).
\end{equation}
If we write~(\ref{hamilton for associated}) using the basis that produced the canonical form~(\ref{canonical form}), the equality~(\ref{fancy q}) guarantees that
\begin{equation}
\label{for later}
\bd^2\widehat{h\slo}(0)=\bd^2h\slo(0)|_{U\slo}=-\mathbb{J}_{4n}\mathcal{A}=
\left(
\begin{array}{cc}
\mathbf{0}&\no\mathbb{J}_{2n}\\
-\no\mathbb{J}_{2n}&-\mathbb{I}_{2n}
\end{array}
\right).
\end{equation}

\medskip

\noindent\textbf{Invariant splitting of the resonance space $U\sno$}Using expressions ~(\ref{hessian of j}) and~(\ref{for later}) we can immediately construct a very convenient splitting of the resonance space $U\sno$:  let $L:U\sno\rightarrow U\sno$ be the linear map defined by 
$\langle L(v), w\rangle=\bd^2(\widehat{h\slo}-\J^{\no})(0)(v, w)$, for any $v, w\in U\sno$. Using expressions~(\ref{hessian of j}),~(\ref{hessian of j matrix}), and~(\ref{for later}) we can write, using the basis introduced in~(\ref{canonical form}), that we will use in all that follows:
\[
L=\left(
\begin{array}{cc}
\mathbf{0}&\mathbf{0}\\
\mathbf{0}&-\mathbb{I}_{2n}
\end{array}
\right).
\]
Since the linear map $L$ is $G\times S^1$--equivariant and self--adjoint  we can split $U\sno=V_0\oplus V_1$ as the direct sum of the two $G\times S^1$--invariant subspaces,
\begin{equation}
\label{v 0 and v 1}
V_0:=\ker L=\left\{\left(\left.\begin{array}{c}
\mathbf{a}\\
\mathbf{0}
\end{array}
\right)\right| \mathbf{a}\in\mathbb{R}^{2n}\right\},\qquad V_1:={\rm Im}\,  L=\left\{\left(\left.\begin{array}{c}
\mathbf{0}\\
\mathbf{b}
\end{array}
\right)
\right| \mathbf{b}\in\mathbb{R}^{2n}\right\}.
\end{equation}
Since by hypothesis~\textbf{(H3)} we are in the generic situation, the resonance space $U\sno$ splits as the sum of two complex dual irreducible subspaces with respect to the $G\times S^1$--representation~\cite{dellnitz melbourne marsden}. Given that by construction $V_0$ and $V_1$ are $G\times S^1$--invariant and have the same dimension, the $G\times S^1$--representations on $V_0$ and $V_1$ are necessarily complex irreducible. We describe more precisely the interplay between the decomposition $U\sno=V_0\oplus V_1$ and the $G\times S^1$--action in the following lemma.

\begin{lemma}
\label{viva euclides!}
In all the matricial statements bellow we assume the use of the basis of the canonical form~(\ref{canonical form}).
\begin{description}
\item[(a)] Let $g\in G\times S^1$ arbitrary and $v=v_0+v_1\in U\sno$, with $v_0\in V_0$ and $v_1\in V_1$. Then, there exists a orthogonal matrix $A_g$ such that $[A_g,\mathbb{J}_{2n}]=0$ and 
\[
g\cdot v=
\left(
\begin{array}{cc}
A_g&\mathbf{0}\\
\mathbf{0}&A_g
\end{array}
\right)\cdot\left(
\begin{array}{c}
v_0\\
v_1
\end{array}
\right).
\]
\item[(b)] The inner product on $U\sno$ that takes the Euclidean form when expressed in the coordinates corresponding to the basis used to write the canonical form~(\ref{canonical form}) is $G\times S^1$--invariant.
\end{description}
\end{lemma}

\noindent\textbf{Proof}\ \textbf{(a)}\  The $G\times S^1$--invariance of the spaces $V_0$ and $V_1$ implies, for any $g\in G\times S^1$, the existence of two invertible matrices $A_g$ and $B_g$ such that
\[
g\cdot v=
\left(
\begin{array}{cc}
A_g&\mathbf{0}\\
\mathbf{0}&B_g
\end{array}
\right)\cdot\left(
\begin{array}{c}
v_0\\
v_1
\end{array}
\right),
\]
for any $v=v_0+v_1\in U\sno$, with $v_0\in V_0$ and $v_1\in V_1$. Given that by hypothesis the family of Hamiltonians $h_\lambda$ is $G$--invariant, the linearization $\mathcal{A}$ in~(\ref{canonical form}) of the vector field $X_{h\slo}$ at the origin is necessarily $G\times S^1$--equivariant and, consequently
\[
\left(
\begin{array}{cc}
A_g&\mathbf{0}\\
\mathbf{0}&B_g
\end{array}
\right)\left(
\begin{array}{cc}
\no \mathbb{J}_{2 n}&\mathbb{I}_{2n}\\
\mathbf{0}&\no \mathbb{J}_{2 n}
\end{array}
\right)=\left(
\begin{array}{cc}
\no \mathbb{J}_{2 n}&\mathbb{I}_{2n}\\
\mathbf{0}&\no \mathbb{J}_{2 n}
\end{array}
\right)\left(
\begin{array}{cc}
A_g&\mathbf{0}\\
\mathbf{0}&B_g
\end{array}
\right),
\]
which implies that $A_g=B_g$ and that $A_g\mathbb{J}_{2n}=\mathbb{J}_{2n}A_g$, for any $g\in G\times S^1$. We now see that the matrices $A_g$ are orthogonal: given that the $G\times S^1$--action is canonical, we have that
\begin{equation}
\label{orthogonal almost}
\left(
\begin{array}{cc}
A_g^T&\mathbf{0}\\
\mathbf{0}&A_g^T
\end{array}
\right)\left(
\begin{array}{cc}
\mathbf{0}&-\mathbb{I}_{2n}\\
\mathbb{I}_{2n}&\mathbf{0}
\end{array}
\right)\left(
\begin{array}{cc}
A_g&\mathbf{0}\\
\mathbf{0}&A_g
\end{array}
\right)=\left(
\begin{array}{cc}
\mathbf{0}&-\mathbb{I}_{2n}\\
\mathbb{I}_{2n}&\mathbf{0}
\end{array}
\right).
\end{equation}
This equality guarantees that $A_g^T A_g=A_g A_g^T=\mathbb{I}_{2n}$, as required.

\medskip

\noindent\textbf{(b)}\  Let $\langle\cdot,\cdot\rangle$ be the inner product in the statement, that is, for any $v, w\in U\sno$ whose coordinates in the basis of~(\ref{canonical form}) are $v=(v_1,\ldots, v_{4n})$ and $w=(w_1,\ldots, w_{4n})$ we have that $\langle v, w\rangle=v_1 w_1+\ldots+v_{4n}w_{4n}$. Notice that, also in this basis, the inner product $\langle\cdot,\cdot\rangle$ can be expressed as
\[
\langle v, w\rangle=\omega|_{U\sno}(\mathbb{J}_{4n}v, w).
\]
Before we proceed notice that expression~(\ref{orthogonal almost}) together with the orthogonality of $A_g$ implies that for any $v\in U\sno$ and any $g\in G$ we have that $\mathbb{J}_{4n}g\cdot v=g\cdot\mathbb{J}_{4n} v$. If we put together this fact with the canonical character of this action we obtain the invariance of the inner product. Indeed, for $v, w\in U\sno$ and $g\in G\times S^1$ arbitrary, we have that
\[
\langle g\cdot v, g\cdot w\rangle=\omega|_{U\sno}(\mathbb{J}_{4n}g\cdot v, g\cdot w)=\omega|_{U\sno}(g\cdot\mathbb{J}_{4n} v, g\cdot w)=\omega|_{U\sno}(\mathbb{J}_{4n} v, w)=\langle v, w\rangle,
\]
which concludes the proof.\ \ \ $\blacksquare$

\begin{remark}
\normalfont
In all our subsequent discussions we will use the inner product presented in the previous lemma and the basis of the canonical form~(\ref{canonical form}).\ \ \ $\blacklozenge$
\end{remark}

\medskip

\noindent\textbf{The quadratic part of the Hamiltonian and a final generic hypothesis} The complex irreducibility of the $G\times S^1$--action on $V_0$ implies~\cite[Lemma 3.4]{g2} that if $\mathcal{P}_{G\times S^1}(V_0)$ denotes the ring of real $G\times S^1$--invariant polynomials on $V_0$, one can choose a basis $\{F_1,\ldots, F_l\}$ of $\mathcal{P}_{G\times S^1}(V_0, V_0)$, that is,  the finite type $\mathcal{P}_{G\times S^1}(V_0)$--module of $G\times S^1$--equivariant polynomial mappings of $V_0$ into itself, such that
\begin{equation}
\label{equivariants}
\left\{
\begin{array}{ll}
F_1=\mathbb{I}_{2n}& \\
F_2=\mathbb{J}_{2n}& \\
{\rm deg}  F_k>1&\forall k>2. 
\end{array}
\right.
\end{equation}
Analogously, one can choose a Hilbert basis $\{\theta_1,\ldots,\theta_r\}$ of the module $\mathcal{P}_{G\times S^1}(V_0)$, such that
\begin{equation}
\label{invariants}
\left\{
\begin{array}{ll}
\theta_1(v)=\|v\|^2&  \\
{\rm deg}  \theta_k>2&\forall k>1. 
\end{array}
\right.
\end{equation}
In particular, the $G\times S^1$--invariance of the Hamiltonians $\widehat{h_\lambda}$ of the equivalent system $(U\sno, \omega|_{U\sno}, \widehat{h_\lambda})$ implies that for each $\lambda$, the second derivative $\bd^2\widehat{h_\lambda}(0)$, considered as a linear map $\bd^2\widehat{h_\lambda}(0):V_0\oplus V_1\rightarrow V_0\oplus V_1$ is $G\times S^1$--equivariant. At the same time, since it is a Hessian, it is symmetric and therefore there are functions $\sigma, \rho, \tau, \psi\in C\suinf(\mathbb{R})$ such that:
\begin{equation}
\label{expression hessian}
\bd^2\widehat{h_\lambda}(0)=\left(
\begin{array}{cc}
\sigma(\lambda)\mathbb{I}_{2n}&\tau(\lambda)\mathbb{I}_{2n}+\psi(\lambda)\mathbb{J}_{2n}\\
\tau(\lambda)\mathbb{I}_{2n}-\psi(\lambda)\mathbb{J}_{2n}&\rho(\lambda)\mathbb{I}_{2n}.
\end{array}\right),
\end{equation}
where, by~(\ref{for later}), we have the following
initial conditions: $\sigma(\lo)=0, \rho(\lo)=-1, \tau(\lo)=0$, and $\psi(\lo)=\no$. In all that follows we will assume the following generic hypothesis:

\begin{description}
\item[(H4)] The one--parameter family of $G$--Hamiltonian systems $(V,\omega,h_\lambda)$ satisfies that $\sigma'(\lo)\neq 0$, where $\sigma(\lambda)\in C\suinf(\mathbb{R})$ is the smooth real function introduced in~(\ref{expression hessian}). 
\end{description}

The previous generic hypothesis that will be of much technical importance in what follows is related, as we show in the following lemma, to the motion of eigenvalues depicted in Figure~\ref{fig:hopf}.

\begin{lemma}
Let $(V,\omega,h_\lambda)$ be a one--parameter family of $G$--Hamiltonian systems satisfying hypotheses~{\rm\textbf{(H1)}} through~{\rm\textbf{(H4)}}. Then, there is a part of the spectrum of the linearization of the Hamiltonian vector fields $X_{h_\lambda}$  at zero that behaves as in Figure~\ref{fig:hopf} as we move the parameter $\lambda$.
\end{lemma}

\noindent\textbf{Proof} We first compute the eigenvalues of  $\mathcal{A}_\lambda$, that is, the restriction  of $A_\lambda$ to the symplectic subspace $(U\sno,\, \omega|_{U\sno})$. This vector field is  Hamiltonian. More specifically: 
\begin{equation}
\label{hamilton for associated restricted}
{\rm{\mathbf i}}_{\mathcal{A}_\lambda}\omega|_{U\sno}=\bd \mathcal{Q}_\lambda,
\end{equation}
where, for any $v\in U\sno$,
\begin{equation}
\label{fancy q restricted}
\mathcal{Q}_\lambda(v)=\frac{1}{2}\bd^2h_\lambda(0)(v,\,v)=\frac{1}{2}\bd^2\widehat{h_\lambda}(0)(v,\,v).
\end{equation}
If we use in~(\ref{fancy q restricted}) the specific form for the symplectic form $\omega$ introduced in~(\ref{canonical form}) we obtain that:
\begin{equation}
\label{relation a and h}
\mathcal{A}_\lambda=\mathbb{J}_{4n}\bd^2\widehat{h_\lambda}(0).
\end{equation}
We now use~(\ref{expression hessian}) and obtain that:
\begin{equation}
\label{a matrix}
\mathcal{A}_\lambda=\left(
\begin{array}{cc}
-\tau(\lambda)\mathbb{I}_{2n}+\psi(\lambda)\mathbb{J}_{2n}&-\rho(\lambda)\mathbb{I}_{2n}\\
\sigma(\lambda)\mathbb{I}_{2n}&\tau(\lambda)\mathbb{I}_{2n}+\psi(\lambda)\mathbb{J}_{2n}
\end{array}\right).
\end{equation}
We compute the eigenvalues of this matrix using the well--known fact~\cite[p. 102, ex. 9]{halmos 74} that if $A, B, C$, and $D$ commute, then
\[
\det\left(
\begin{array}{cc}
A&B\\
C&D
\end{array}
\right)=\det(AD-BC).
\]
Indeed, using this relation, we see that the characteristic polynomial of $\mathcal{A}_\lambda$ is
\begin{align*}
\det(\mathcal{A}_\lambda-&\mu\mathbb{I}_{4n})=\left(
\begin{array}{cc}
(-\tau(\lambda)-\mu)\mathbb{I}_{2n}+\psi(\lambda)\mathbb{J}_{2n}&-\rho(\lambda)\mathbb{I}_{2n}\\
\sigma(\lambda)\mathbb{I}_{2n}&(\tau(\lambda)-\mu)\mathbb{I}_{2n}+\psi(\lambda)\mathbb{J}_{2n}
\end{array}\right)\\
	&=\det\left[(\mu^2-\tau(\lambda)^2-\psi(\lambda)^2+\rho(\lambda)\sigma(\lambda))\mathbb{I}_{2n}-2\mu\psi(\lambda)\mathbb{J}_{2n}\right]\\
	&=\det\left(
\begin{array}{cc}
(\mu^2-\tau(\lambda)^2-\psi(\lambda)^2+\rho(\lambda)\sigma(\lambda))\mathbb{I}_n&2\mu\psi(\lambda)\mathbb{I}_n\\
-2\mu\psi(\lambda)\mathbb{I}_n&(\mu^2-\tau(\lambda)^2-\psi(\lambda)^2+\rho(\lambda)\sigma(\lambda))\mathbb{I}_n
\end{array}\right)\\
	&=\det\left[((\mu^2-\tau(\lambda)^2-\psi(\lambda)^2+\rho(\lambda)\sigma(\lambda))^2+(2\mu\psi(\lambda))^2)\mathbb{I}_n\right]\\
	&=\left(((\mu^2-\tau(\lambda)^2-\psi(\lambda)^2+\rho(\lambda)\sigma(\lambda))^2+4\mu^2\psi(\lambda)^2\right)^n.
\end{align*}
Consequently, the eigenvalues of $\mathcal{A}_\lambda$ are
\begin{equation}
\label{eigenvalues of a at last}
\mu(\lambda)=\pm\sqrt{\tau(\lambda)^2-\rho(\lambda)\sigma(\lambda)-\psi(\lambda)^2\pm 2|\psi(\lambda)|\sqrt{\rho(\lambda)\sigma(\lambda)-\tau(\lambda)^2}}.
\end{equation}

We are now in position to show that hypothesis~\textbf{(H4)}
implies the motion depicted in Figure~\ref{fig:hopf}. Indeed, assume that $\sigma'(\lo)\neq 0$. Let $f_1(\lambda)=\rho(\lambda)\sigma(\lambda)-\tau(\lambda)^2$. Note that $f_1(\lo)=0$ and $f_1'(\lo)=-\sigma'(\lo)$, that by hypothesis is different from zero. This implies that the function $f_1$ changes sign at $\lo$. More explicitely, suppose  that $\sigma'(\lo)<0$ (the case $\sigma'(\lo)>0$ is completely analogous); in that case $f_1(\lambda)<0$ for $\lambda<\lo$ and $f_1(\lambda)>0$ for $\lambda>\lo$. Since we have that:
\[\mu(\lambda)=\pm\sqrt{-f_1(\lambda)-\psi(\lambda)^2\pm 2|\psi(\lambda)|\sqrt{f_1(\lambda)}}=\pm\sqrt{-(\sqrt{f_1(\lambda)}\pm |\psi(\lambda)|)^2}.\]
then, by considering the cases $\lambda<\lo$, $\lambda=\lo$, and $\lambda>\lo$ in the previous expression taking into account the changes of sign in $f_1$, we obtain the evolution of eigenvalues  illustrated in  Figure~\ref{fig:hopf}.\ \ \ $\blacksquare$

\begin{remark}
\normalfont
The generic hypothesis \textbf{(H4)} is a sufficient but not necessary condition for obtaining a behavior of the eigenvalues as in Figure~\ref{fig:hopf}, that is, such an evolution can take place even for systems in which $\sigma'(\lo)=0$.\ \ \ $\blacklozenge$ 
\end{remark}

\medskip

\section{Hamiltonian Hopf bifurcation and relative periodic orbits}
\label{the main result}

The main goal of this section is the statement and proof of a result that will provide and estimate on the number of relative periodic orbits of a one--parameter family of $G$--Hamiltonian systems $(V,\omega,h_\lambda)$ that satisfies the hypotheses~\textbf{(H1)} through~\textbf{(H4)}, formulated in the previous section.

We will begin by introducing some classical definitions that will make more explicit some of the concepts used in the previous paragraphs.

It appears very frequently in examples dealing with symmetric families of Hamiltonian systems  that the canonical symmetry group $G$ contains a continuous {\bfi globally Hamiltonian symmetry}: suppose that $G$ contains a Lie subgroup $H$ of positive dimension. We say that the canonical action of $H$ on $V$ is globally Hamiltonian when we can associate to it an equivariant momentum map $\K:V\rightarrow\hs$ which is defined by the fact that its components $\K\suxi:=\langle\K,\xi\rangle\in C\suinf(\mathbb{R})$, $\xi\in\h$, have as associated Hamiltonian vector fields the infinitesimal generators of the action
\[
\xi_V(v)=\ddto\exp t\xi\cdot v\qquad\xi\in\h, v\in V.
\]

\begin{definition}
\label{rpo}
Let $(V, \omega, h)$ be a Hamiltonian system with a symmetry given by the canonical action of the Lie group $H$ on $V$. The point 
$v\in V$ is called a {\bfi relative periodic point (RPP)},
if there is a $\tau>0$ and an element $g\in H$ such that
\[
F_{t+\tau}(v)=g\cdot F_t(v)\qquad\text{for any}\qquad
t\in\mathbb{R},
\] 
where $F_t$ is the flow of the Hamiltonian vector
field $X_h$. The set
\[
\gamma(v):=\{F_t(v)\ \mid\ t>0\}
\]
is called a {\bfi relative periodic orbit (RPO)} through $v$. 
The constant $\tau>0$ is its {\bfi relative period\/} and the
group element $g\in H$ is its {\bfi phase shift\/}.
\end{definition}

\begin{proposition}
\label{rpo and j}
Let $(V, \omega, h)$ be a Hamiltonian system with a globally Hamiltonian symmetry given by the canonical action of the Lie group $H$ on $V$ with associated momentum map $\K:V\rightarrow\hs$. If the Hamiltonian vector field $X_{h-\K\suxi}$, $\xi\in\h$, has a periodic point $v\in V$ with period $\tau$, then the point $v$ is a RPP of $X_h$ with relative period $\tau$ and phase shift $\exp \tau\xi$.
\end{proposition}

\noindent\textbf{Proof} Let $F_t$ be the flow of the Hamiltonian vector field $X_h$ and $K_t(v)=\exp t\xi\cdot v$ that of $X_{\K\suxi}$. By Noether's Theorem:
\[
[X_h,X_{\K\suxi}]=-X_{\{h,\K\suxi\}}=0,
\]
where the bracket $\{\cdot,\cdot\}$ denotes the Poisson bracket associated to the symplectic form $\omega$. Due to this equality, we can write (see for instance~\cite[Corollary 4.1.27]{mta}) the following expression for $G_t$, the flow of $X_{h-\K\suxi}$:
\[
G_t(v)=\lim_{n\rightarrow\infty}(F_{t/n}\circ K_{-t/n})^n(v)=(K_{-t}\circ F_t)(v)=\exp -t\xi\cdot F_t(v).
\]
Since by hypothesis the point $v$ is periodic for $G_t$ with period $\tau$, we have that
\[
v=\exp -\tau\xi\cdot F_\tau(v),
\]
or, equivalently,
\[
F_{\tau}(v)=\exp \tau\xi\cdot v,
\]
as required.\ \ \ $\blacksquare$

\medskip

Using the previous proposition, we will reduce the search for RPOs of a generic one--parameter family of $G$--Hamiltonian systems $(V, \omega, h_\lambda)$ that satisfies conditions~{\rm\textbf{(H1)}}, {\rm\textbf{(H2)}},~{\rm\textbf{(H3)}}, and~{\rm\textbf{(H4)}}, to the search for periodic orbits of the vector fields of the form $X_{h_\lambda-\K\suxi}$, and will prove the following result:

\begin{theorem}
\label{relative periodic theorem}
Let $(V, \omega, h_\lambda)$ be a  one--parameter family of $G$--Hamiltonian systems that satisfies conditions~{\rm\textbf{(H1)}}, {\rm\textbf{(H2)}},~{\rm\textbf{(H3)}}, and~{\rm\textbf{(H4)}}. Suppose that $G$ contains a Lie subgroup $H$ of positive dimension with associated equivariant momentum map $\K:V\rightarrow\hs$. Let  $U\sno$ be the resonance space with primitive period $T\sno$. Then, for each $\xi\in\h$ whose norm $\|\xi\|$ is small enough, there are at least, in each energy level  nearby zero and for each value of the parameter $\lambda$ near $\lo$, as many
relative periodic orbits as the number of equilibria of a $G\suxi\times S^1$--equivariant vector field defined on the unit sphere on $V_0$. The symbol  $G\suxi$ denotes the adjoint isotropy subgroup of the element $\xi\in\h$, that is,
\[
G\suxi=\{g\in G\mid \Ad_g\xi=\xi\}.
\]
\end{theorem}

\begin{remark}
\normalfont
If we are just interested in looking for purely periodic orbits it suffices to use  Theorem~\ref{relative periodic theorem} with  $\xi=0$. Conversely, if we use this result with a value of the parameter $\xi\neq 0$ we cannot conclude that the predicted RPOs are not trivial, that is, that they are not just periodic orbits. This point will become much clearer in the  examples presented in the following sections.\ \ \ $\blacklozenge$ 
\end{remark}

\begin{remark}
\normalfont
In terms of practical applications, the relevance of Theorem~\ref{relative periodic theorem} is given by the fact that the estimate that it provides in terms of the number of equilibria of an equivariant vector field on the sphere can sometimes be calculated via topological arguments, as we will see later on.\ \ \ $\blacklozenge$ 
\end{remark}

\noindent\textbf{Proof} We will work in the basis of the resonance space $U\sno$ provided by the equivariant Williamson normal form, in particular we will use the matricial expressions~(\ref{canonical form}), which are consistent with the decomposition $U\sno=V_0\oplus V_1$ presented in~(\ref{v 0 and v 1}). Recall that the subspaces $V_0$ and $V_1$ are $G\times S^1$--invariant. Abusing the notation a little bit we will use the symbol $\xi$ to denote both an element of the Lie algebra $\h\subset\g$ and its representation on $V_0$ and $V_1$. Using Lemma~\ref{viva euclides!} we can write, for each $v=v_0+v_1\in U\sno$ represented in the previously mentioned basis, 
\[
\xi_{U\sno}(v)=\xi\cdot v_0+\xi\cdot v_1=
\left(
\begin{array}{rr}
\xi&\mathbf{0}\\
\mathbf{0}&\xi
\end{array}
\right)
\left(
\begin{array}{c}
v_0\\
v_1
\end{array}
\right).
\]
Note that, also  by Lemma~\ref{viva euclides!},  that the matrix $\xi$ is skew--symmetric, $\xi^T=-\xi$, therefore normal, and hence diagonalizable. The same Lemma implies that the linear map $\xi:V_0\rightarrow V_0$ associated to $\xi\in\h$ commutes with $\mathbb{J}_{2n}$, $[\xi,\mathbb{J}_{2n}]=0$, and consequently these two endomorphisms can be simultaneously diagonalized.

We recall that,
\[
\langle\K(v),\xi\rangle=\frac{1}{2}\omega(\xi\cdot v,v)=\frac{1}{2}(v_0,v_1)
\left(
\begin{array}{rr}
\mathbf{0}&\xi\\
-\xi&\mathbf{0}
\end{array}
\right)
\left(
\begin{array}{c}
v_0\\
v_1
\end{array}
\right).
\]
In particular,
\[
\bd^2\K\suxi(0)=\left(
\begin{array}{rr}
\mathbf{0}&\xi\\
-\xi&\mathbf{0}
\end{array}
\right).
\]

We start the proof by defining the $\mathbb{R}\times\h$--parameter family of Hamiltonian functions given by
\[
h_{\lambda,\xi}=h_\lambda-\K\suxi.
\]
Due to the hypotheses on the family $h_\lambda$, the quadratic nature of the momentum map $\K$, and the fact that $h_{\lambda,0}=h_\lambda$, the family $h_{\lambda,\xi}$ satisfies the hypotheses of the Normal Form Reduction Theorem~\cite{vanderbauwhede 95}. Therefore, a new family $\widehat{h_{\lambda,\xi}}$ can be constructed such that, for any value $(\lambda,\xi)$ of the parameters, the Hamiltonian $\widehat{h_{\lambda,\xi}}$ is $S^1$--invariant with respect to the action generated by the semisimple part of the linearization at zero of 
$X_{h_{\lo,0}}=X_{h\slo}$, that is, $(\theta,v)\mapsto{\rm e}^{\frac{\theta}{\no}\mathcal{A}^s} v$, $\theta\in S^1$, with 
\begin{equation}
\label{matrix of s 1}
\mathcal{A}^s=
\left(
\begin{array}{cc}
\no\mathbb{J}_{2n}&\mathbf{0}\\
\mathbf{0}&\no\mathbb{J}_{2n}
\end{array}
\right).
\end{equation}
The Normal Form Reduction Theorem guarantees that the $S^1$--relative equilibria of $\widehat{h_{\lambda,\xi}}$ are in correspondence with the periodic orbits $h_{\lambda,\xi}$ which, by Proposition~\ref{rpo and j}, are RPOs of $h_\lambda$. The quadratic nature of the momentum map $\K$ and its $S^1$--invariance imply that $\widehat{h_{\lambda,\xi}}$ can be chosen to be of the form 
\[
\widehat{h_{\lambda,\xi}}=\widehat{h_\lambda}-\K\suxi,
\]
with $\widehat{h_\lambda}$ the normal form for the family $h_\lambda$. 

As a result of these premises, the RPOs that we are looking for will be given by  the critical points of the function $\widehat{h_\lambda}-\K\suxi-\J^{\zeta+\alpha}$, that is, the elements $(v,\alpha,\lambda,\xi)\in U\sno\times\mathbb{R}\times\mathbb{R}\times\h$ for which the function
\begin{equation}
\label{equation in brute}
F^\zeta(v,\alpha,\lambda):=\nabla_{U\sno}\left(\widehat{h_\lambda}-\K\suxi-\J^{\zeta+\alpha}\right)(v)
\end{equation}
has a zero. As customary, the gradient in the previous expression is constructed using the inner product introduced in Lemma~\ref{viva euclides!}.

\medskip

\noindent\textbf{Lyapunov--Schmidt reduction and the bifurcation equation} The linearization $L^\zeta:U\sno\rightarrow U\sno$ of the equation~(\ref{equation in brute}) at the point $(0,0,\lo,0)$ produces, in the usual basis, the expression:
\begin{equation}
\label{lyapunov or not}
L^\zeta=\bd^2\left(\widehat{h_\lambda}-\J^{\zeta}\right)(0)=\bd^2\left(h_\lambda-\J^{\zeta}\right)(0)=\left(
\begin{array}{cc}
\mathbf{0}&(1-\frac{\zeta}{\no})\no\mathbb{J}_{2n}\\
-(1-\frac{\zeta}{\no})\no\mathbb{J}_{2n}&-\mathbb{I}_{2n}
\end{array}
\right).
\end{equation}
By looking at this matricial expression we see that it is possible to Lyapunov--Schmidt reduce the bifurcation problem posed in~(\ref{equation in brute}) whenever  $\zeta=\no$, which we will assume in the sequel. In those circumstances $\ker L^{\no}=V_0$, ${\rm Im}\,   L^{\no}=V_1$. Let $\mathbb{P}:U\sno\rightarrow V_0$ be the $G\times S^1$--equivariant projection associated to the splitting $U\sno=V_0\oplus V_1$. The equation $(\mathbb{I}-\mathbb{P})F^{\no}(v_0+v_1, \alpha, \lambda,\xi)=(\mathbb{I}-\mathbb{P})\nabla_{U\sno}\left(\widehat{h_\lambda}-\K\suxi-\J^{\no+\alpha}\right)(v_0+v_1)=0$ defines,
via the Implicit Function Theorem, a  function $v_1:V_0\times\mathbb{R}\times \mathbb{R}\times\h\rightarrow V_1$, such that
\begin{equation}
\label{v 1 r p o}
(\mathbb{I}-\mathbb{P})F^{\no}(v_0+v_1(v_0, \alpha, \lambda,\xi), \alpha, \lambda,\xi)=(\mathbb{I}-\mathbb{P})\nabla_{U\sno}\left(\widehat{h_\lambda}-\K\suxi-\J^{\no+\alpha}\right)(v_0+v_1(v_0, \alpha, \lambda,\xi))=0.
\end{equation}
Notice that the function $v_1(v_0, \alpha, \lambda,\xi)$ is $G\suxi\times S^1$--equivariant, since this is the symmetry under which $F^{\no}$ is equivariant, that is, for any $g\in G\suxi\times S^1$, we have that $v_1(g\cdot v_0, \alpha, \lambda,\xi)=g\cdot v_1(v_0, \alpha, \lambda,\xi)$.

The final Lyapunov--Schmidt $G\suxi\times S^1$--equivariant {\bfi reduced bifurcation equation}, whose zeros provide us with the RPOs that we are after, is given by $B:V_0\times\mathbb{R}\times \mathbb{R}\times\h\rightarrow V_0$, where
\begin{eqnarray}
B(v_0,\alpha,\lambda,\xi)&=&\mathbb{P}F^{\no}(v_0+v_1(v_0,\alpha,\lambda,\xi),\alpha,\lambda,\xi)=\mathbb{P}\nabla_{U\sno}\left(\widehat{h_\lambda}-\K\suxi-\J^{\no+\alpha}\right)(v_0+v_1(v_0,\alpha,\lambda,\xi))\notag\\
	&=&\nabla_{U\sno}\left(\widehat{h_\lambda}-\K\suxi-\J^{\no+\alpha}\right)(v_0+v_1(v_0, \alpha, \lambda,\xi))\qquad\text{(by~(\ref{v 1 r p o})).}\label{bifurcation equation r p o}
\end{eqnarray}

We collect the main properties of the reduced bifurcation equation in the following

\begin{lemma}
\label{gradient}
The reduced bifurcation equation~{\rm (\ref{bifurcation equation r p o})} is $G\suxi\times S^1$--equivariant with respect to the action of this Lie group on $V_0$ and it is the gradient of a $G\suxi\times S^1$--invariant function defined on $V_0$, that is,
\[
B(v_0, \alpha, \lambda,\xi)=\nabla_{V_0}g(v_0, \alpha, \lambda,\xi),
\] 
where the function $g:V_0\times{\rm Lie}(S^1)\times \mathbb{R}\times\h\rightarrow V_0$  is defined by
\[
g(v_0, \alpha, \lambda,\xi)=(\widehat{h_\lambda}-\J^{\no+\alpha}-\K\suxi)(v_0+v_1(v_0, \alpha, \lambda,\xi)).
\]
\end{lemma}

\noindent\textbf{Proof} The $G\suxi\times S^1$--equivariance is a direct consequence of the construction of $B$. As to the gradient character of $B$, note first that for any $w\in V_1$ we have that
\begin{align}
\langle F^{\no}(v_0+v_1(v_0, \alpha, \lambda,\xi), \alpha,\lambda,\xi), w\rangle&=\langle F^{\no}(v_0+v_1(v_0, \alpha,\lambda,\xi), \alpha,\lambda,\xi), (\mathbb{I}-\mathbb{P})w\rangle\notag\\
	&=\langle (\mathbb{I}-\mathbb{P})F^{\no}(v_0+v_1(v_0, \alpha,\lambda,\xi), \alpha,\lambda,\xi), w\rangle=0\label{intermediate}
\end{align}
where the last equality follows from the construction of the function $v_1$ through expression~(\ref{v 1 r p o}). Now, let $u\in V_0$ arbitrary. We write:
\begin{align*}
\langle B(v_0, \alpha,\lambda,\xi), u\rangle&=\langle \mathbb{P}F^{\no}(v_0+v_1(v_0, \alpha,\lambda,\xi), \alpha,\lambda,\xi), u\rangle\\
	&=\langle F^{\no}(v_0+v_1(v_0, \alpha,\lambda,\xi), \alpha,\lambda,\xi), u\rangle\\
	&=\langle F^{\no}(v_0+v_1(v_0, \alpha,\lambda,\xi), \alpha,\lambda,\xi), u+D_{V_0}v_1(v_0, \alpha,\lambda,\xi)\cdot u\rangle\quad(\text{by }~(\ref{intermediate}))\\
	&=\langle \nabla_{U\sno}(\widehat{h_\lambda}-\J^{\no+\alpha}-\K\suxi)(v_0+v_1(v_0, \alpha,\lambda,\xi)), u+D_{V_0}v_1(v_0, \alpha,\lambda,\xi)\cdot u\rangle\\
	&=\bd(\widehat{h_\lambda}-\J^{\no+\alpha}-\K\suxi)(v_0+v_1(v_0, \alpha,\lambda,\xi))\cdot (u+D_{V_0}v_1(v_0, \alpha,\lambda,\xi)\cdot u)\\
	&=\bd g(v_0, \alpha,\lambda,\xi)\cdot u=\langle\nabla_{V_0}g(v_0, \alpha,\lambda,\xi), u \rangle,
\end{align*}
as required. 
This construction is a particular case of the one carried out in~\cite{constrained paper} and~\cite{pascal}.\ \ \ $\blacktriangledown$

\medskip

\noindent\textbf{Notational simplification:} In order to make notation a little bit lighter we will assume in the rest of the proof, without loss of generality,  that the system has been scaled in such a way that $\no=1$ and $\lo=0$. 

\medskip

\noindent The following lemmas provide a local description of the reduced bifurcation equation that will be much needed.

\begin{lemma}
\label{derivative v 1}
The function $v_1$ introduced in~{\rm (\ref{v 1 r p o})} has the following  two properties:
\begin{align}
\text{{\rm \textbf{(i)}}}\ & v_1(0,\alpha,\lambda,\xi)=0\qquad\text{for all}\qquad \alpha, \lambda\in\mathbb{R}, \text{\ and}\ \xi\in\h.\\
\text{{\rm \textbf{(ii)}}}\ & D_{V_0}v_1(0,\alpha,\lambda,\xi)=-\frac{\tau(\lambda)}{\rho(\lambda)}\mathbb{I}_{2n}-\frac{(1+\alpha)-\psi(\lambda)}{\rho(\lambda)}\mathbb{J}_{2n}-\frac{1}{\rho(\lambda)}\xi.\label{v 1 r p o o}
\end{align}
\end{lemma}

\noindent\textbf{Proof} Part~\textbf{(i)} is a consequence of the uniqueness of the solutions provided by the Implicit Function Theorem. The proof of part~\textbf{(ii)} is supplied in the Appendix, Section~\ref{lemma 1}.\ \ \ $\blacktriangledown$

\medskip

\noindent The proof of the following lemma is a lengthy but straightforward computation.
\begin{lemma}
\label{principal part b}
Let $B(v_0,\alpha,\lambda,\xi)$ be the reduced bifurcation equation, then:
\begin{description}
\item[(i)] 
\[
\!\!\!\!\!\!\!\!\!\!\!\!\!\!\!D_{V_0}B(0,\alpha,\lambda,\xi)=
\frac{\sigma(\lambda)\rho(\lambda)-\tau^2(\lambda)-((1+\alpha)-\psi(\lambda))^2}{\rho(\lambda)}\mathbb{I}_{2n}+\frac{2\left[(1+\alpha)-\psi(\lambda)\right]}{\rho(\lambda)}\mathbb{J}_{2n}\xi+\frac{\xi^2}{\rho(\lambda)}.
\]
\item[(ii)] The principal part of the reduced bifurcation equation is given by the expression:
\begin{multline}
\label{principal part bifurcation}
B(v_0,\alpha,\lambda,\xi)=(\lambda\sigma'(0)+\alpha^2)v_0-\xi^2 v_0-2\alpha\mathbb{J}_{2n}\xi v_0-2\psi'(0)\alpha\lambda v_0\\
+2\psi'(0)\lambda\mathbb{J}_{2n}\xi v_0+C\left(v_0^{(3)}\right)+\text{h.o.t.},
\end{multline}
where $C\left(v_0^{(3)}\right)$ is the trilinear operator obtained by taking the gradient of the fourth order term in the $v_0$--expansion of $\widehat{h\slo}(v_0+v_1(v_0,0,0,0))$.
\end{description}
\end{lemma}

We now write the reduced bifurcation equation in polar coordinates, that is, we define
\[
B_p(r,u_0,\alpha,\lambda,\xi)=B(ru_0,\alpha,\lambda,\xi),
\]
where $r\in\mathbb{R}$ and $u_0\in S^{\dim V_0-1}$. We introduce the function
\[
F(r,u_0,\alpha,\lambda,\xi)=\frac{\langle B(ru_0,\alpha,\lambda,\xi),u_0\rangle}{r}.
\]
By looking at~(\ref{principal part bifurcation}) it is clear that the function $F$ is smooth at the origin, $F(0,0,0,0,0)=0$ and that $D_\lambda  F(0,0,0,0,0)=\sigma'(0)\neq 0$, by hypothesis~\textbf{(H4)}. Therefore, the Implicit Function Theorem guarantees the existence of a smooth function $\lambda(r,u_0,\alpha,\xi)$ such that $\lambda(0,0,0,0)=0$ and $F(r,u_0,\alpha,\lambda(r,u_0,\alpha,\xi),\xi)=0$. This equality implies that if we substitute the function $\lambda(r,u_0,\alpha,\xi)$ on the reduced bifurcation equation, this time considered as a vector field on $V_0$, we obtain a new $(\alpha,\xi)$--parameter dependent vector field
\begin{equation}
\label{g sphere}
G(r,u_0,\alpha,\xi)=B_b(r,u_0,\alpha,\lambda(r,u_0,\alpha,\xi),\xi)
\end{equation}
which due to the fact that $\langle B_b(r,u_0,\alpha,\lambda(r,u_0,\alpha,\xi),\xi),u_0\rangle=0$ is, for each small enough fixed value of $r$, a $G\suxi\times S^1$--equivariant vector field on the sphere on $V_0$ of radius $r$, whose zeroes constitute solutions of the reduced bifurcation equation. 
\ \ \ $\blacksquare$

\medskip

\medskip

\noindent\textbf{Method for the optimal use of Theorem~\ref{relative periodic theorem}}
The optimal and most organized way to apply Theorem~\ref{relative periodic theorem} consists of using the estimate it provides in the fixed point subspaces $V_0^H$ corresponding to the various  subgroups  $H$ in the lattice of isotropy subgroups of the $G\times S^1$--action on $V_0$, replacing the group $G\times S^1$ by  $N(H)$, which is a group that acts on $V_0^{H}$ (not necessarily in an irreducible manner). The symbol $N(H)$ denotes the normalizer of $H$ in $G\times S^1$ and $V_0^{H}$ is the vector subspace of $V_0$ formed by the vectors fixed by $H$.  We make more explicit this comment in the following paragraphs.

Let $H$ be a subgroup  of $G\times S^1$. If $\pi:G\times S^1\rightarrow G$ denotes the canonical projection and $\pi(H)=:K\subset G$, Proposition 7.2 in~\cite{g2} guarantees the existence of a group homomorphism $\theta:K\rightarrow S^1$ such that 
\begin{equation}
\label{character subgroup}
H=\{(k,\theta(k))\in G\times S^1\mid k\in K\}.
\end{equation}
In our discussion we will be concerned with {\bfi spatiotemporal symmetries}, that is, subgroups $H$ of $G\times S^1$ for which the homomorphism $\theta:K\rightarrow S^1$ is nontrivial. Using the characterization~(\ref{character subgroup}) it is straightforward to see that
\[
N(H)=N_G(K)\times S^1.
\]
The $N_G(K)$--action on $U\sno^H$ is globally Hamiltonian with momentum map $\K^H:U\sno^H\rightarrow {\rm Lie}\left(N_G(K)\right)\sus$ given by the restriction of the $G$--momentum map to $U\sno^H$, that is, for any $v\in U\sno^H$ and any $\xi\in {\rm Lie}\left(N_G(K)\right)$, we have that
\[
\langle \K^H(v),\xi\rangle=\langle \K(v),\xi\rangle.
\]
The same statement applies to the $S^1$--action. Using these objects we can reformulate Theorem~\ref{relative periodic theorem}  on the fixed point spaces $V_0^H$.

\begin{corollary}
\label{non maximal remark}
Let $(V, \omega, h_\lambda)$ be a  one--parameter family of $G$--Hamiltonian systems that satisfies conditions~{\rm\textbf{(H1)}}, {\rm\textbf{(H2)}},~{\rm\textbf{(H3)}}, and~{\rm\textbf{(H4)}}. Let $H$ be a spatiotemporal isotropy subgroup  of the $G\times S^1$--action on $V_0$, such that $\dim V_0^H=2k$ and $K:=\pi(H)$. Then, for each $\xi\in {\rm Lie}\left(N_G(K)\right)$ whose norm $\|\xi\|$ is small enough, there are at least
in each energy level  nearby zero and for each value of the parameter $\lambda$ near $\lo$, as many
relative periodic orbits as the number of equilibria of a $N_G(K)\suxi\times S^1$--equivariant vector field on the unit sphere on $V_0^H$.
The relative periods of these RPOs are close to $T\sno$, and their phase shifts are close to $\exp T\sno\xi$. The symbol $N_G(K)\suxi$ denotes the adjoint isotropy subgroup of the element $\xi\in{\rm Lie}\left(N_G(K)\right)$, that is,
\[
N_G(K)\suxi=\{g\in N_G(K)\mid \Ad_g\xi=\xi\}.
\]
The mapping $\pi:G\times S^1\rightarrow G$ denotes the canonical projection.
\end{corollary}

\medskip

As we already said, both the previous result and Theorem~\ref{relative periodic theorem} can be used to look for purely periodic motions by taking in their respective statements $\xi=0$. There is a situation of special interest:

\medskip

\noindent\textbf{Periodic orbits with maximal isotropy subgroup:} Let $H$ be a maximal isotropy subgroup of the $G\times S^1$--action on $V_0$. In the presence of maximality we have at our disposal the following convenient result:
\begin{lemma}
\label{maximal advantages}
Let $H$ be a maximal isotropy subgroup of the compact $G\times S^1$--action on $V_0$. Let $N$ be the Lie group $\ele$ and $N^0$ be the connected component of the identity of $N$. Then either
\begin{description}
\item[(i)] $N^0\simeq S^1$, and $N/N^0=\{Id\}$ or $N/N^0\simeq\mathbb{Z}_2$, or
\item[(ii)] $N^0\simeq SU(2)$ and $N\simeq SU(2)$.
\end{description}
In the first case we say that $H$ is a maximal complex subgroup. In the second case we say that $H$ a maximal quaternionic subgroup.
\end{lemma}

\noindent\textbf{Proof} It is a straightforward combination of the general result for  linear actions of compact Lie groups~\cite{bredon, golubitsky 83, g2} with Proposition 12.5 in~\cite{golubitsky stewart hopf} that eliminates the possibility of having real maximal isotropy subgroups when the compact group in question is $G\times S^1$.\ \ \ $\blacksquare$

\medskip

Using the previous lemma and an additional genericity hypothesis, the estimate given in Theorem~\ref{relative periodic theorem} can be made very explicit:

\begin{corollary}
\label{hamiltonian fiedler}
Let $(V,\,\omega,\,h_\lambda)$ be a generic one--parameter family of $G$--Hamiltonian systems that satisfies conditions~{\rm\textbf{(H1)}}, {\rm\textbf{(H2)}}, {\rm\textbf{(H3)}}, and~{\rm\textbf{(H4)}}. Let $H$ be a maximal isotropy subgroup  of the $G\times S^1$--action on $V_0$ such that $\dim (V_0^H)=l\neq 0$. Then:
\begin{description}
\item[(i)] If $N^0\simeq S^1$ there are at least $l/2$ (if $N/N^0=\{Id\}$) or $l/4$ (if $N/N^0\simeq\mathbb{Z}_2$) branches of periodic solutions with isotropy $H$ coming out of the origin as one varies the parameter $\lambda$, with periods close to $T\sno$. 
\item[(ii)] If $N^0\simeq SU(2)$ there are at least $l/4$  branches of periodic solutions with isotropy $H$ coming out of the origin as one varies the parameter $\lambda$, with periods close to $T\sno$.
\end{description}
\end{corollary}

\noindent\textbf{Proof}  We will adapt to our problem the approach followed in~\cite{chossat maximal, muriel} for rotating waves. The main idea behind the proof consists of using the maximality hypothesis to give a numerical evaluation of the estimate in Corollary~\ref{non maximal remark}, that is, the number of equilibria of a $N(H)/H$--equivariant vector field on the sphere $S^{l-1}$.

More especifically, let
$G^H:=G|_{V_0^H}$ be the restriction of the vector field $G$ on $V_0$, defined in~(\ref{g sphere}), to the  fixed point set $V_0^H$, and $G_r^H(u_0,\alpha):=G^H(r,u_0,\alpha)$ be the $N$--equivariant vector field on $S_r^{l-1}$ obtained by fixing $r$ in the mapping $G^H$ (note that in our case $\xi=0$ since we are looking for periodic orbits). The zeroes of this vector field are in one to one correspondence with the solutions that we search. Due to the maximality hypothesis on the subgroup $H$, the $N$--action on the sphere $S_r^{l-1}$ is free and therefore the corresponding orbit space $S_r^{l-1}/N$ is a smooth manifold onto which we can project the $N$--equivariant vector field $G_r^H$. Let $\bar{G}_r^H$ be the projected vector field. Due to the genericity hypothesis in the statement, the Poincar\'e--Hopf Theorem allows us to say that $\bar{G}_r^H$ has at least $\chi(S_r^{l-1}/N)$ equilibria, where $\chi$ denotes the Euler characteristic. These zeroes lift to equilibria of the restriction of the reduced bifurcation equation to $V_0^H$, due to the gradient character (see Lemma~\ref{gradient}) of  $B$ and consequently of its restriction to $V_0^H$.

In order to conclude our argument it is enough to show that $\chi(S^{l-1}/N)$ corresponds to the estimates provided in the statement of the theorem. In the first case, when $N^0\simeq S^1$, the dimension of $V_0^H$ is necessarily even (we will write $l=2k$ for certain $k\in\mathbb{N}$) and there are two possibilities: the quotient $N/N^0$ is either $\{Id\}$ or it is isomorphic to $\mathbb{Z}_2$. If $N/N^0=\{Id\}$:
\[\chi(S^{l-1}/N)=\chi((S^{l-1}/N^0)/(N/N^0))=\chi(S^{2k-1}/S^1)=\chi(\mathbb{CP}^{k-1})=k=\frac{l}{2}.\]
If $N/N^0\simeq\mathbb{Z}_2$:
\[\chi(S^{l-1}/N)=\chi((S^{l-1}/N^0)/(N/N^0))=\chi(\mathbb{CP}^{k-1}/\mathbb{Z}_2)=\frac{k}{2}=\frac{l}{4},
\]
where we used the well--known fact that if $G$ is a finite group acting freely on a manifold $M$, then (see for instance~\cite[Corollary 5.22]{kawakubo})
\[
\chi\left(\frac{M}{G}\right)=\frac{\chi{(M)}}{|G|}.
\]
Finally, if $H$ is maximal quaternionic then $l=\dim(V_0^H)=4k$ for some $k\in\mathbb{N}$, necessarily, and
\[\chi(S^{l-1}/N)=\chi(S^{4k-1}/SU(2))=\chi(\mathbb{HP}^{k-1})=k=\frac{l}{4}.
\]
The calculation of the Euler characteristic $\chi(\mathbb{HP}^{k-1})$ of the quaternionic projective space is made using an argument based the spectral series of Leray (see for instance~\cite{bott}). 

The computations that we just carried out give us periodic orbits for a fixed $r$. Moving smoothly this parameter we obtain the branches required in the statement of the theorem. \ \ \ $\blacksquare$

\section{Bifurcation of non-periodic relative periodic orbits in the presence of extra hypotheses}

The tools presented in Theorem~\ref{relative periodic theorem} for the search of RPOs based on topological methods produce estimates that, as we will see in the following examples, have some limitations, in particular we have no example where it guarantees the bifurcation of non periodic RPOs. This circumstance has motivated us to use a more analytical approach under dimensional hypotheses that are satisfied in very relevant situations. A detailed study of the bifurcation equation in the presence of these hypotheses will provide us with sharper estimates that completely describe all the bifurcation phenomena that we see in the examples.

\subsection{Motivating example: two coupled harmonic oscillators subjected to a magnetic field as an example of symmetric Hamiltonian Hopf bifurcation}
\label{example}

We consider the system formed by two identical particles with unit charge in the plane, subjected to identical harmonic forces, to a homogeneous magnetic field perpendicular in direction to the plane of motion, and to an interaction potential that will preserve certain group of symmetry. We will denote by $(q_1,q_2)$ the coordinates of the configuration space of the first particle and by $(q_3,q_4)$ those of the second one. If $\gamma$ is a constant that determines the intensity of the magnetic field, it is easy to see that the Hamiltonian function of the system described above is
\begin{multline}
\label{hamiltonian harmonic}
H(\mathbf{q},\mathbf{p})=\frac{1}{2m}(p_1^2+p_2^2+p_3^2+p_4^2)+\left(\frac{\gamma^2}{2m}-\frac{k}{2}\right)(q_1^2+q_2^2+q_3^2+q_4^2)\\
+\frac{\gamma}{m}(p_1q_2-p_2q_1)+\frac{\gamma}{m}(p_3q_4-p_4q_3)+f(\pi_1^i,\pi_2^i,\pi_3^i),
\end{multline}
where
\[
\pi_1^i=q_i^2+q_{i+2}^2,\quad\pi_2^i=p_i^2+p_{i+2}^2,\quad\pi_3^i=p_iq_{i+2}-p_{i+2}q_i,\quad\pi_4^i=q_i p_i+q_{i+2}p_{i+2},\quad i\in\{1,2\},
\]
and $f$ is a higher order function on its variables that expresses a non linear interaction between the two particles.

This system has, for all values of the parameters $\gamma$ and $k$, an equilibrium at the point $(q_1,q_3,q_2,q_4,p_1,p_3,p_2,p_4)=(\mathbf{0},\mathbf{0})$. The linearization of the dynamics at that point is represented by the matrix (the coordinates are ordered as in the previous equality)
\begin{equation}
\label{linear harmonic magnetic}
\mathcal{A}_k=\left(
\begin{array}{cc}
-\frac{\gamma}{m}\mathbb{J}_4&\frac{1}{m}\mathbb{I}_4\\
\left(k-\frac{\gamma^2}{m}\right)\mathbb{I}_4&-\frac{\gamma}{m}\mathbb{J}_4
\end{array}
\right),
\end{equation}
whose eigenvalues are
\[
\lambda_k=\pm\frac{1}{m}\sqrt{km-2\gamma^2\pm 2\gamma\sqrt{\gamma^2-km}}.
\]
If we move the parameter $k$ around the value $k_\circ=\gamma^2/m$ these eigenvalues present a Hamiltonian Hopf behavior like the one depicted in Figure~\ref{fig:hopf}. 

We now study the symmetries of the system. Note that after the assumptions on the interaction function $f$, the system is invariant under the canonical $S^1$--action given by the lifted action to the phase space of 
\[
(\varphi,\mathbf{q})\longmapsto
\left(
\begin{array}{cc}
\begin{array}{cc}
\cos\varphi&-\sin\varphi\\
\sin\varphi&\cos\varphi
\end{array}&\mbox{\protect\LARGE 0}\\
\mbox{\protect\LARGE 0}&\begin{array}{cc}
\cos\varphi&-\sin\varphi\\
\sin\varphi&\cos\varphi
\end{array}
\end{array}
\right)\cdot\mathbf{q},
\]
where $\mathbf{q}=(q_1,q_3,q_2,q_4)$, and by the transformation  
\[
\tau\cdot
\left(
\begin{array}{c}
q_1\\
q_2\\
q_3\\
q_4
\end{array}
\right)=\left(
\begin{array}{c}
q_1\\
q_2\\
-q_3\\
-q_4
\end{array}
\right).
\]
The momentum map $\K:\mathbb{R}^8\rightarrow\mathbb{R}$ associated to the $S^1$--action is given by the expression $\K(\mathbf{q},\mathbf{p})=p_3q_1-q_3p_1-p_2q_4+p_4q_2$.

If we now look at the linearization~(\ref{linear harmonic magnetic}) evaluated at the Hopf value $k_\circ=\gamma^2/m$ we see that in this case $V_0$ consists of the points of the form $(q_1,q_3,q_2,q_4,\mathbf{0})$. The $S^1$--action on $V_0$ generated by the semisimple part of  $\mathcal{A}_{k_\circ}$ can be written as
\[
(\theta,\mathbf{q})\longmapsto {\rm e}^{-\theta\frac{\gamma}{m}\mathbb{J}_4}\cdot\mathbf{q}.
\]
In order to better study the  group actions on $V_0$ we will perform a linear change of variables. Let $(z_1,z_2,\bar{z_1},\bar{z_2})$ be the new (complex) coordinates, given by
\begin{equation}
\label{complex change}
\begin{array}{rl}
z_1&=q_1+q_4+iq_2-iq_3\\
z_2&=q_1-q_4+iq_2+iq_3\\
\bar{z_1}&=q_1+q_4-iq_2+iq_3\\
\bar{z_2}&=q_1-q_4-iq_2-iq_3.
\end{array}
\end{equation}
If we take as new angles $\psi_1$ and $\psi_2$, defined by:
\[
\psi_1=\varphi+\frac{\gamma}{m}\theta,\qquad \psi_2=\varphi-\frac{\gamma}{m}\theta,
\]
we realize that the previously introduced actions form a $O(2)\times S^1$--action on $V_0$ that takes the following convenient simple expression:
\begin{equation}
\label{perfect action}
(\psi_1,\psi_2)\cdot(z_1,z_2)=({\rm e}^{i\psi_1}z_1,{\rm e}^{i\psi_2}z_2)\qquad\text{and}\qquad \tau\cdot(z_1,z_2)=(z_2,z_1).
\end{equation}
That is, we have shown that the system of two coupled harmonic oscillators subjected to a magnetic field, whose Hamiltonian is given by~(\ref{hamiltonian harmonic}) can be taken as an example of Hamiltonian Hopf bifurcation with $O(2)\times S^1$--symmetry, which we will study in full generality in the following subsection.

\subsection{RPOs in Hamiltonian Hopf bifurcation with $O(2)$--symmetry}

Having the example in the previous section as a motivation we will study in what follows the RPOs that appear in a Hamiltonian Hopf bifurcation phenomenon in the presence of a $O(2)$--symmetry. The simplicity of this symmetry will allow us to explicitly write down the principal part of the reduced bifurcation equation in full generality, and to read off directly from it the RPOs that we are looking for.

We start by recalling that in the canonical coordinates introduced in~(\ref{canonical form}) the principal part of the reduced bifurcation equation is, by Lemma~\ref{principal part b}, equal in our case to
\begin{multline*}
B(v_0,\alpha,\lambda,\xi)=(\lambda\sigma'(\lo)+\alpha^2\no^2)v_0-\xi^2 v_0-2\alpha\no\mathbb{J}_{4}\xi v_0-2\psi'(\lo)\no\alpha\lambda v_0\\
+2\psi'(\lo)\lambda\mathbb{J}_{4}\xi v_0+\mathbb{P}\bd^4 h_{\lo}(0)\left(v_0^{(4)}\right)+\text{h.o.t.}.
\end{multline*}
We now rewrite this expression in the coordinates $(z_1,z_2,\bar{z_1},\bar{z_2})$ in which the $O(2)\times S^1$--action looks like~(\ref{perfect action}). In doing so we need to express in these new coordinates the matrix $\mathbb{J}_{4}\xi$, which can be easily achieved by using the explicit expression of the change of variables~(\ref{complex change}). Indeed, we have that in those coordinates
\[
\mathbb{J}_{4}\xi\equiv
\left(
\begin{array}{cccc}
1&0&0&0\\
0&-1&0&0\\
0&0&1&0\\
0&0&0&-1
\end{array}
\right),
\]
and therefore, the first terms of the expansion of the reduced bifurcation equation are:
\begin{equation}
\label{bifurcation o 2}
B(\mathbf{z},\alpha,\lambda,\xi)=(\lambda\sigma'(\lo)+\alpha^2\no^2)\left(\begin{array}{c}
z_1\\
z_2\\
\bar{z_1}\\
\bar{z_2}
\end{array}
\right)-2\no\alpha\xi\left(\begin{array}{r}
z_1\\
-z_2\\
\bar{z_1}\\
-\bar{z_2}
\end{array}
\right)+\left(\begin{array}{c}
(a|z_1|^2+b|z_2|^2)z_1\\
(a|z_2|^2+b|z_1|^2)z_2\\
(a|z_1|^2+b|z_2|^2)\bar{z_1}\\
(a|z_2|^2+b|z_1|^2)\bar{z_2}
\end{array}
\right)+\ldots,
\end{equation}
where the coefficients $a$ and $b$ are related to the fourth order terms in the expansion of the Hamiltonian, that is, $\mathbb{P}\bd^4 h_{\lo}(0)\left(v_0^{(4)}\right)$. In order to keep the simplicity of the exposition we will assume that these two coefficients are non zero and non equal (otherwise we would have to go to higher orders in expression~(\ref{bifurcation o 2})). The RPOs that we are looking for are given by the solutions of the system of equations:
\begin{align}
0&=(\lambda\sigma'(\lo)+\alpha^2\no^2) z_1-2\no\alpha\xi z_1+(a|z_1|^2+b|z_2|^2)z_1+\ldots\label{o 2 1}\\
0&=(\lambda\sigma'(\lo)+\alpha^2\no^2)z_2+2\no\alpha\xi z_2+(a|z_2|^2+b|z_1|^2)z_2+\ldots\label{o 2 2}
\end{align}
Since $\sigma'(\lo)\neq 0$, equation~(\ref{o 2 1}) can be easily solved by dividing the expression by $z_1$ and then using the Implicit Function Theorem to define a function 
\begin{equation}
\label{lambda o 2}
\lambda\equiv\lambda(z_1,z_2,\alpha,\xi)=\frac{1}{\sigma'(\lo)}\left[-\alpha^2\no^2+2\no\alpha\xi-(a|z_1|^2+b|z_2|^2)+\ldots\right]
\end{equation}
that substituted into~(\ref{o 2 1}) solves it. Hence, plugging~(\ref{lambda o 2}) into~(\ref{o 2 2}) we reduce the problem to solving a scalar equation which can be done again via the Implicit Function Theorem: we divide~(\ref{o 2 1}) by $z_1$ and~(\ref{o 2 2}) by $z_2$ and obtain: 
\begin{equation}
\label{we will see}
0=4\no\alpha\xi+(a-b)(|z_2|^2-|z_1|^2)+\ldots
\end{equation}
Since for equivariance reasons $z_1$ and $z_2$ always appear in the tail of the previous expression as combinations of $|z_1|^2$ and $|z_2|^2$, the hypotheses on the coefficients $a$ and $b$ allow us to solve this final scalar equation by defining a function
\begin{equation}
\label{final z}
|z_2|^2\equiv|z_2|^2\left(|z_1|^2,\alpha,\xi\right)=|z_1|^2-\frac{4\no\alpha\xi}{a-b}+\ldots
\end{equation}
\begin{figure}[htb]
\begin{center}
\includegraphics{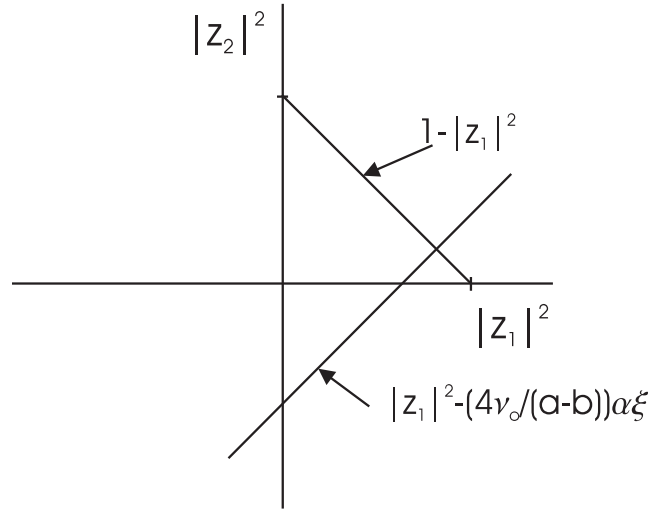}
\caption{Parameterization of the RPOs in the Hamiltonian Hopf bifurcation with $O(2)$--symmetry.}\label{fig:fixed norm}
\end{center}
\end{figure}
As we illustrate in Figure~\ref{fig:fixed norm}, the solution~(\ref{final z}) predicts, for each fixed value of the norm $|z_1|^2+|z_2|^2$ a one--parameter family of RPOs that are obtained by making vary the product $\alpha\xi$. More explicitly, and using Figure~\ref{fig:fixed norm}, suppose that $|z_1|^2+|z_2|^2=1$ and that the value $\alpha\xi$ is fixed, then, the intersection of the lines $|z_2|^2=1-|z_1|^2$ and $|z_2|^2=|z_1|^2-\frac{4\no\alpha\xi}{a-b}$ provides us with the abovementioned RPO.

Notice that all these RPOs cannot be predicted by merely using Theorem~\ref{relative periodic theorem} since even though all the hypotheses needed in this result are fulfilled, it only predicts two RPOs for each value of the norm $|z_1|^2+|z_2|^2$. Moreover, one cannot be sure that these are not just periodic motions since the predicted orbits could lie in fixed spaces of maximal isotropy, thereby implying their periodicity.

\begin{remark}
\normalfont
In contrast with the Hamiltonian case, the Hopf bifurcation of nontrivial RPOs in the dissipative case with $O(2)$ symmetry is subjected to the presence of a cubic order degeneracy in the normal form. Unfolding this singularity leads to a codimension two bifurcation problem where the RPOs appear as a secondary branching from the primary branches of periodic orbits.\ \ \ $\blacklozenge$
\end{remark}
 
\medskip

\subsection{Hamiltonian Hopf bifurcation of RPOs for reduced integrable systems}
\label{manual treatment}

The analysis performed in the previous section dealing with the Hamiltonian Hopf bifurcation with $O(2)$--symmetry is a particular case of a more general situation. Indeed, the main feature in that example that allowed us to carry out a by hand in--depth study of the reduced bifurcation equation was the coincidence of one half the dimension of the reduced space $V_0$ with the dimension of the symmetry group $O(2)\times S^1$. We will see in this section that whenever we are in the presence of this \emph{reduced integrability hypothesis} an analysis in the same style can be performed.

More explicitly, all along  this section we will be dealing with $(V, \omega, h_\lambda)$, a  one--parameter family of $G$--Hamiltonian systems that satisfies conditions~{\rm\textbf{(H1)}}, {\rm\textbf{(H2)}},~{\rm\textbf{(H3)}}, and~{\rm\textbf{(H4)}} such that if  $4n$ is the dimension of the resonance space $U\sno$ with primitive period $T\sno$, then the rank of $G\times S^1$ equals $n$, that is, the maximal tori of the Lie group $G\times S^1$ have all dimension equal to $n$. 

Let $\mathbb{T}^{n-1}\subset G$ be a maximal torus of $G$, and let  $\xi\in\frak{t}^{n-1}$ be an element in the Lie algebra of $\mathbb{T}^{n-1}$. 
Like in the previous section we can find coordinates in which the action of $\mathbb{T}^{n-1}\times S^1$ looks simple. Namely, there exists a set of complex coordinates $(z_1,\ldots,z_n,\bar{z_1},\ldots,\bar{z_n})$ (and conjugates) for $V_0$ and a set of angular coordinates $(\xi_1,\ldots,\xi_{n-1})$ for the torus $T^{n-1}$, for which the  $\mathbb{T}^{n-1}$--action looks like
\[
({\rm e}^{i\xi_1},\ldots,{\rm e}^{i\xi_{n-1}})\cdot(z_1,\ldots,z_n)=({\rm e}^{i\xi_1}z_1,\ldots,{\rm e}^{i\xi_{n-1}}z_{n-1},{\rm e}^{i(c_1\xi_1+\ldots+c_{n-1}\xi_{n-1})}z_n),
\]
where the coefficients $c_1,\ldots,c_{n-1}$ are rational constants. If we incorporate the  $S^1$--action using these complex coordinates, the $\mathbb{T}^{n-1}\times S^1$--action looks like
\[
({\rm e}^{i xi_1},\ldots,{\rm e}^{i \xi_{n-1}},{\rm e}^{i \alpha})\cdot(z_1,\ldots,z_n)=({\rm e}^{i(\xi_1+\alpha)}z_1,\ldots,{\rm e}^{i(\xi_{n-1}+\alpha)}z_{n-1},{\rm e}^{i(c_1\xi_1+\ldots+c_{n-1}\xi_{n-1}+\alpha)}z_n).
\]
Let us now set
\begin{eqnarray}
\psi_j &=& \xi_j + \alpha~,~~j=1,\dots,n-1 \\
\psi_n &=& c_1\xi_1+\cdots +c_{n-1}\xi_{n-1}+\alpha.
\end{eqnarray}
Under the condition
\begin{equation} \label{condition}
c_1 + \cdots +c_{n-1} \neq 1
\end{equation}
these relations define a change of coordinates on the $n$-dimensional torus $\mathbb{T}^{n-1}\times S^1$, and in these new coordinates the action can now be written in the very simple fashion
\begin{equation} \label{action_canonique}
({\rm e}^{i\psi_1},\ldots,{\rm e}^{i\psi_{n}})\cdot(z_1,\ldots,z_n)=({\rm e}^{i\psi_1}z_1,\ldots,{\rm e}^{i\psi_n}z_n).
\end{equation}
Notice that under Condition (\ref{condition}) the ring of invariant polynomials for this action on $V_0$ is generated by the quadratic invariants $\pi_j=z_j\bar z_j$, $j=1,\dots,n$, and that the strata of this action are obtained by setting some of the $z_j$'s equal to 0 while keeping the others different from 0. The orbit space for this action can be identified with the positive cone in $\mathbb{R}^n$ $\{(\pi_1,\dots,\pi_n)~/~\pi_j\geq 0,~j=1,\dots,n\}$.

Recall now that in the canonical coordinates introduced in~(\ref{canonical form}), the principal part of the reduced bifurcation equation is, by Lemma~\ref{principal part b}, equal to
\begin{multline}
\label{thing to transform}
B(v_0,\alpha,\lambda,\xi)=(\lambda\sigma'(\lo)+\alpha^2\no^2)v_0-\xi^2 v_0-2\alpha\no\mathbb{J}_{2n}\xi v_0-2\psi'(\lo)\no\alpha\lambda v_0\\
+2\psi'(\lo)\lambda\mathbb{J}_{2n}\xi v_0+C\left(v_0^{(3)}\right)+\text{h.o.t.},
\end{multline}

From (\ref{action_canonique}) it is clear that the matrices $\mathbb{J}_{2n}$, $\mathbb{J}_{2n}\xi $ and $\xi^2$ in~(\ref{thing to transform}) take, in these newly introduced coordinates, the form:
\[
\mathbb{J}_{2n}=
\left(
\begin{array}{ccc}
i&\cdots&0\\
\vdots&\ddots&\vdots\\
0&\cdots&i
\end{array}
\right),\quad
\mathbb{J}_{2n}\xi=
\left(
\begin{array}{cccc}
-\psi_1&\cdots&\cdots&0\\
\vdots&\ddots& &\vdots\\
\vdots& &-\psi_{n-1}&\vdots\\
0&\cdots&\cdots&-\psi_n
\end{array}
\right),
\]
\[
\xi^2=
-\left(
\begin{array}{cccc}
\psi_1^2&\cdots&\cdots&0\\
\vdots&\ddots& &\vdots\\
\vdots& &\psi_{n-1}^2&\vdots\\
0&\cdots&\cdots&\psi_n^2
\end{array}
\right).
\]
Using these new coordinates and the $\mathbb{T}^{n-1}\times S^1$ equivariance of $B$, we rewrite the components of~(\ref{thing to transform})  (we omit the complex conjugate part). For $i\in\{1,\ldots,n\}$ we have:
$$
B_i(z,\alpha,\lambda,\xi)=\left(\lambda\sigma'(\lo)+\psi_i^2 + \widehat{C_i}\left(|z_1|^2,\ldots,|z_n|^2\right)  + \text{h.o.t.}\right) z_i
$$
where 
$$
\widehat{C_i}\left(|z_1|^2,\ldots,|z_n|^2\right)=\widehat{c_{i1}}|z_1|^2+\ldots+\widehat{c_{in}}|z_n|^2,\qquad \widehat{c_{ij}}\in\mathbb{R},\quad i,j\in\{1,\ldots,n\},
$$
We can now state a theorem about the bifurcation of RPOs.

\begin{theorem} \label{theorem_rpos}
Let $(V, \omega, h_\lambda)$ be a one--parameter family of $G$--Hamiltonian systems that satisfies conditions~{\rm\textbf{(H1)}}, {\rm\textbf{(H2)}},~{\rm\textbf{(H3)}}, and~{\rm\textbf{(H4)}}.
Suppose that: (i) the dimension of $V_0$ equals twice the rank $n$ of $G\times S^1$; (ii) the condition (\ref{condition}) on the torus action is satisfied. Then, if the matrix
$$
\Delta = \left( \widehat{c_{nj}}-\widehat{c_{ij}} \right),~~1\leq i\leq n,~1\leq j\leq n-1,
$$
has maximal rank $n-1$, there exists a family of RPOs with $n$ different frequencies which bifurcates from the trivial solution as $\lambda$ crosses $\lambda_0$.
\end{theorem}
\begin{remark}
\normalfont
The condition on the matrix $\Delta$ is not generic, because the values of the coefficients $\widehat{c_{ij}}$ are constrained by the $G$--equivariance of the operator $C$.
\end{remark}

\noindent\textbf{Proof}  Since we are looking for solutions with $z_i\neq 0$, for all $i$, we can factor out $z_i$ in each equation $B_i=0$. The resulting equations read, for $i\in\{1,\dots,n\}$,
$$
0 = \lambda\sigma'(\lo)+\psi_i^2 + \widehat{C_i}\left(|z_1|^2,\ldots,|z_n|^2\right) + \text{h.o.t.}
$$
These equations are simply those that we would have obtained by projecting first (\ref{thing to transform}) on the orbit space corresponding to the toral action; in all that follows we will denote $|z_i|^2$ by $\pi_i$.
Since by hypothesis \textbf{(H4)}, $\sigma'(\lo)\neq 0$, we can solve any one of these equations for $\lambda$. Let us do so for the equation with $i=n$. By substituting $\lambda$ by the resulting expression in the remaining equations, we have reduced the problem to solving a system of $n-1$ equations which, at leading order, have the form
\begin{equation} \label{eq_orbitspace}
0 = \psi_i^2 - \psi_n^2 + (\widehat{c_{i1}}-\widehat{c_{n1}})\pi_1 + \cdots + (\widehat{c_{in}}-\widehat{c_{nn}})\pi_n + \text{h.o.t.}.
\end{equation} 
If the matrix $\Delta$ has maximal rank, we obtain a unique family of solutions of the system (\ref{eq_orbitspace}), for which $n-1$ of the $\pi_i$'s depend smoothly on the remaining one and on the parameters $\psi_j$, $j=1,\dots,n$. In order to fix notations, let us assume without loss of generality that we have obtained $\pi_i=\pi_i(\psi_1,\dots,\psi_n,\pi_n)$ for $i=1,\dots,n-1$. These solutions still have to lie inside the orbit space, that is,  we still have to check the additional conditions $\pi_i\geq 0$. However, since the $\psi_j$'s are free parameters, the quantities $\psi_i^2-\psi_n^2$ can take any real value. Therefore, if we set $z_n=0$, we can always find values for the $\psi_j$'s such that $\pi_i>0$ for $i=1,\dots,n-1$. These inequalities are still satisfied if the $\psi_j$'s are close enough to these values and $\pi_n>0$ is close enough to 0. Finally, since these solutions lie on the principal stratum for the action of $\mathbb{T}^{n-1}\times S^1$, the corresponding RPOs have $n$ different frequencies which depend smoothly on $\lambda$.  $\blacksquare$

\begin{remark}
\normalfont
In the problem with $O(2)$ symmetry analyzed in Section 4.2, we have $n=2$, $c_1=-1$, $\widehat{c_{11}}=\widehat{c_{22}}$ and $\widehat{c_{12}}=\widehat{c_{21}}$. The hypotheses of  Theorem~\ref{theorem_rpos} are therefore generically satisfied in this case.\ \ \ $\blacklozenge$
\end{remark}
\begin{remark}
\normalfont
As it was already the case with Theorem~\ref{relative periodic theorem}, Theorem \ref{theorem_rpos} still applies if instead of $G\times S^1$ acting in $V_0$, we consider the group $N(H)/H$ acting in $V_0^H$, where $H$ is an isotropy subgroup of the $G\times S^1$--action. In the next section we shall see an application of this remark.\ \ \ $\blacklozenge$
\end{remark}

\subsection{Hamiltonian Hopf bifurcation with $SO(3)$ symmetry}
\label{spheric example}

Hopf bifurcation problems with $SO(3)$ symmetry for dissipative systems have been investigated by several authors in the case in which the eigenspaces associated with the critical eigenvalues is the direct sum of twice the five dimensional (real) irreducible representation of $SO(3)$ (see \cite{golubitsky stewart hopf}), \cite{IoRo}, \cite{mrs}, and \cite{Leis}). This is the simplest possible case with $SO(3)$ symmetry which does not reduce to Hopf bifurcation with either trivial or $O(2)$ symmetry. Nevertheless, it leads to a normal form in a ten dimensional real vector space. The list of solutions with maximal isotropy, hence purely periodic ones, has been given in \cite{golubitsky stewart hopf} and in \cite{mrs}. However, the most interesting feature of this problem is the possibility of having a bifurcated branch of RPOs in a six dimensional subspace. This was first found by \cite{IoRo}. Another approach was taken by \cite{Leis} (using orbit space reduction) who did not recover the result of \cite{IoRo}. This remark shows the level of difficulty found in obtaining direct branching of RPOs via Hopf bifurcation for equivariant vector fields. In the Hamiltonian context, a related work by Haaf, Roberts and Stewart \cite{Haaf-Roberts-Stewart} has shown the existence of families of periodic orbits with maximal isotropy for a Hamiltonian in $\mathbb{R}^{10}$ which is invariant under the same $SO(3)$--action. 

In the sequel we investigate the Hamiltonian Hopf bifurcation with $SO(3)$ symmetry, when the subspaces $V_0$ and $V_1$ are associated with this ten dimensional representation and we shall see that, in this case, Theorem~\ref{theorem_rpos} applies and shows the existence of several families of RPOs. 

Let 
$$V=V_0\oplus V_1\simeq \mathbb{R}^{10}.$$
We identify $\mathbb{R}^{10}$ with $\mathbb{R}^{5}\otimes\mathbb{C}$ and consider the action of $SO(3)$ on $\mathbb{R}^{5}$ given by its irreducible representation on the space of spherical harmonics of degree 2. Equivalently, we may identify $\mathbb{R}^{5}$ with the space $W$ of $3\times 3$ real symmetric matrices with trace $0$, and consider the action of $SO(3)$ on $W$ defined by
$$\rho_A(M)=A^{-1}MA,~~A\in SO(3), ~M\in W.$$
This definition extends naturally to $\mathbb{R}^{5}\otimes\mathbb{C}$ with the same formula, $M$ now having complex coefficients. We shall therefore identify in all that follows $V_0$ with $W\otimes\mathbb{C}$. The $S^1$--action on $V_0$ is simply defined as multiplication by $e^{i\theta}$ in $\mathbb{C}$, that is, $\theta\cdot M := e^{i\theta}M$.

Any $M\in W\otimes\mathbb{C}$ decomposes uniquely as
$$M =\sum_{m=-2}^2 z_mB_m$$
where
\begin{eqnarray}
B_0 &=& \left(\begin{array}{ccc} 1 & 0 & 0 \\ 0 & 1 & 0 \\ 0 & 0 & 2 \end{array}\right),~~
B_1 = \left(\begin{array}{ccc} 0 & 0 & 1 \\ 0 & 0 & i \\ 1 & i & 0 \end{array}\right),~~B_{-1}=\overline{B}_1 \\
B_{2} &=& \left(\begin{array}{ccc} 1 & i & 0 \\ i & -1 & 0 \\ 0 & 0 & 0 \end{array}\right),~~B_{-2}=\overline{B}_2.
\end{eqnarray}

\begin{figure}[htb]
\begin{center}
\includegraphics{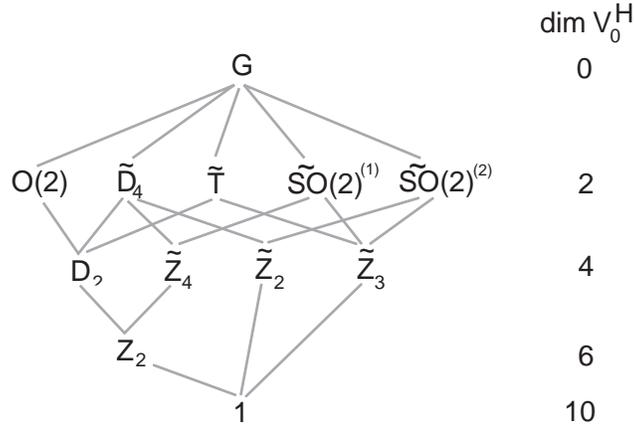}
\caption{Isotropy lattice of the $SO(3)\times S^1$--action.}\label{lattice-O3}
\end{center}
\end{figure}

We now list the isotropy types for the action of $G=SO(3)\times S^1$ on $V_0$ that we have just defined.  We use the presentation and results of \cite{Haaf-Roberts-Stewart}. Figure~\ref{lattice-O3} shows the isotropy lattice of the $G$--action and the dimension of the corresponding fixed--point subspaces.
\\
{\bf Notations:} $\tilde H$ is a subgroup isomorphic to $H\subset O(3)$ but such that $\tilde H \cap S^1 \neq 1$ (here $1$ is the trivial group). In particular, $\tilde\mathbb{Z}_n$ is the group generated by $(R_n,-2\pi/n)\in SO(3)\times S^1$, where $R_n$ is a rotation of angle $2\pi/n$. 

By Corollary~\ref{non maximal remark}, for any subgroup $H$ with $dim(V_0^H)=2$ there exists a branch of periodic solutions having this symmetry. Let us now consider those isotropy subgroups having a fixed-point subspace of dimension 4. The largest subgroup acting faithfully in $V_0^H$ is $N(H)/H$. We list below the different $N(H)/H$ for the non maximal isotropy subgroups:
\begin{enumerate}
\item $N(D_2)/D_2 \simeq D_3\times S^1$;
\item $N(\tilde\mathbb{Z}_4)/\tilde\mathbb{Z}_4 \simeq N(\tilde\mathbb{Z}_2)/\tilde\mathbb{Z}_2 \simeq O(2)\times S^1$;
\item $N(\tilde\mathbb{Z}_3)/\tilde\mathbb{Z}_3 \simeq SO(2)\times S^1$;
\item $N(\mathbb{Z}_2)/\mathbb{Z}_2\simeq O(2)\times S^1$;
\item $N(1)=G$.
\end{enumerate}
In Case 1 we see that solutions with isotropy $D_2$ are always periodic. Case 2 corresponds to the problem described in Section 4.2 (Hopf bifurcation with $O(2)$ symmetry). It was noticed in \cite{Haaf-Roberts-Stewart} that the equations in $V_0^H$ do not degenerate despite the fact that they come from a system with higher symmetry, which leads us to conclude that families of RPOs with two frequencies and with spatio-temporal symmetry $\tilde\mathbb{Z}_4$ as well as $\tilde\mathbb{Z}_2$ do generically bifurcate.
Case 3  falls in the framework of Section 4.3, since the symmetry group is $SO(2)\times S^1$. However, because the equations in $V_0^{\tilde\mathbb{Z}_3}$ are the restriction in that subspace of a system with higher symmetry in $V_0$, we need to compute the cubic order terms in order to insure that no "hidden" degeneracy occurs. We use the argument proved in~\cite{Haaf-Roberts-Stewart}; we can choose $\tilde\mathbb{Z}_3$ so that, by introducing complex coordinates,
\begin{equation}
V_0^{\tilde\mathbb{Z}_3} = \{z_1B_1 + z_2B_{-2},~(z_1,z_2)\in\mathbb{C}^2\}
\end{equation}
and the action of $SO(2)$ is then defined by
$$\phi\cdot (w,z)=(e^{i\phi},e^{-2i\phi}).$$
With the notations of Section 4.3, we therefore have 
$$\psi_1=\phi~\mbox{ and }~c_1=-2.$$
The expression for the cubic $G$--equivariant terms is
$$C(M^{(3)}) = b_1tr(M\bar M)M + b_2tr(M^2)\bar M + b_3\left(M^2\bar M+\bar MM^2-\frac{2}{3}tr(M^2\bar M)Id\right);$$
with $b_j$ real coefficients depending on the specific Hamiltonian at hand. 
Setting $M=z_1B_1+z_2B_{-2}$ we obtain after calculation in $V_0^{\tilde\mathbb{Z}_3}$ the expression
$$C(M^{(3)}) = 4\left((b_1+\frac{b_3}{2})|z_1|^2+(b_1+b_3)|z_2|^2\right)z_1B_1 + 4\left((b_1+b_3)|z_1|^2+b_1|z_2|^2\right)z_2B_{-2}.$$
Let us now check whether the $1\times 2$ matrix $\Delta$ of Theorem \ref{theorem_rpos} has maximal rank. From the above expression we deduce that
$$
\widehat{c_{11}}-\widehat{c_{21}}=-2b_3,~~\widehat{c_{12}}-\widehat{c_{22}}=4b_3.
$$
Therefore the maximality hypothesis is satisfied iff $b_3\neq 0$ (which is a generic condition). \\
Cases 4 and 5 are beyond the range of applicability of Theorem \ref{theorem_rpos}.

\section{Appendix}
\label{appendix}
\subsection{On the invariance properties of the resonance subspace}

In what follows we will sketch the proof of some of the facts about the invariance properties of the resonance subspace mentioned in the preliminaries section when $(V, \omega)$ is a symplectic 
representation space of the Lie group $G$ 
and the Hamiltonian vector field $A$ is 
$G$--equivariant.

\medskip

\noindent\textbf{The resonance subspace $U\sno$ is $G$--invariant} Let $A=A_s+A_n$ be the Jordan--Chevalley decomposition of $A$. Since by hypothesis $A$ is $G$--equivariant, if $\Phi:G\times V\rightarrow V$ denotes the $G$--action, for any $g\in G$, we have that $\Phi_gA=A\Phi_g$. Equivalently, $\Phi_gA_s+\Phi_gA_n=A_s\Phi_g+A_n\Phi_g$, and hence $\Phi_gA_n\Phi_{g\inv}+\Phi_gA_s\Phi_{g\inv}=A_n+A_s$. Since $\Phi_gA_n\Phi_{g\inv}$ is nilpotent, $\Phi_gA_s\Phi_{g\inv}$ is semisimple, $[\Phi_gA_n\Phi_{g\inv},  \Phi_gA_s\Phi_{g\inv}]=0$, and the Jordan--Chevalley decomposition is unique, we have that \[\Phi_gA_n\Phi_{g\inv}=A_n\qquad \text{and}\qquad \Phi_gA_s\Phi_{g\inv}=A_s,\]
necessarily. This implies the $G$--invariance of $U\sno=\ker(e^{A_sT\sno}-I)$. Indeed, let $v\in U\sno$. Hence,  $e^{A_sT\sno}v=v$. At the same time, for any $g\in G$, \[e^{A_sT\sno}(\Phi_gv)=\Phi_ge^{A_sT\sno}v=\Phi_gv,\] 
hence $\Phi_gv\in U\sno$, that is, $U\sno$ is $G$--invariant.

\medskip

\noindent\textbf{The $S^1$--action and the $G$--action on $U\sno$ commute} Let $\Psi:S^1\times U\sno\rightarrow U\sno$ be the $S^1$--action on $U\sno$. For any $g\in G$ and any $\theta\in S^1$:
\[\Phi_g\Psi_\theta=\Phi_ge^{\theta A_s}=e^{\theta A_s}\Phi_g=\Psi_\theta\Phi_g,\]
as required.

\subsection{Proof of Lemma~\ref{derivative v 1}}
\label{lemma 1}
The defining relation~(\ref{v 1 r p o}) of the function $v_1$ implies that for any $v_0\in V_0$, $\alpha, \lambda\in\mathbb{R}$,  and $\xi\in\h$ we have that 
\[
(\mathbb{I}-\mathbb{P})\nabla_{U\sno}(\widehat{h_\lambda}-\J^{1+\alpha}-\K\suxi)(v_0+v_1(v_0, \alpha, \lambda,\xi))=0.
\]
Consequently, for any $w_1\in V_1$:
\begin{eqnarray*}
0&=&\ddto\langle\nabla_{U\sno}(\widehat{h_\lambda}-\J^{1+\alpha}-\K\suxi)(tv_0+v_1(tv_0, \alpha, \lambda,\xi)), w_1\rangle\\
	&=&\ddto \bd (\widehat{h_\lambda}-\J^{1+\alpha}-\K\suxi)(tv_0+v_1(tv_0, \alpha, \lambda,\xi))\cdot w_1\\
	&=&\bd^2(\widehat{h_\lambda}-\J^{1+\alpha}-\K\suxi)(0)(v_0+D_{V_0}v_1(0, \alpha, \lambda,\xi)\cdot v_0, w_1).
\end{eqnarray*}
If we use~(\ref{hessian of j matrix}) and~(\ref{expression hessian}), the previous expression can be matricially expressed as
\begin{eqnarray*}
0&=&(0, w_1)\cdot\left(
\begin{array}{cc}
\sigma(\lambda)\mathbb{I}_{2n}&\tau(\lambda)\mathbb{I}_{2n}+(\psi(\lambda)-(1+\alpha))\mathbb{J}_{2n}-\xi\\
\tau(\lambda)\mathbb{I}_{2n}-(\psi(\lambda)-(1+\alpha))\mathbb{J}_{2n}+\xi&\rho(\lambda)\mathbb{I}_{2n}
\end{array}
\right)\\
	& &\cdot\left(\begin{array}{cc}
v_0\\
D_{V_0}v_1(0, \alpha, \lambda,\xi)\cdot v_0
\end{array}\right)\\
	&=&w_1^T\left[\tau(\lambda)\mathbb{I}_{2n}-(\psi(\lambda)-(1+\alpha))\mathbb{J}_{2n}+\xi+\rho(\lambda)D_{V_0}v_1(0, \alpha, \lambda,\xi)\right]v_0.
\end{eqnarray*}
Given that the previous equation is valid for no matter what $v_0\in V_0$ and $w_1\in V_1$, we can conclude that
\[
D_{V_0}v_1(0, \alpha, \lambda,\xi)=\frac{\psi(\lambda)-(1+\alpha)}{\rho(\lambda)}\mathbb{J}_{2n}-\frac{\tau(\lambda)}{\rho(\lambda)}\mathbb{I}_{2n}-\frac{\xi}{\rho(\lambda)},
\]
as required.\ \ \ $\blacksquare$

\bigskip

\noindent\textbf{Acknowledgments.}  We thank M. Dellnitz and I. Melbourne for their help and patience concerning our questions on their equivariant Williamson normal form~\cite{melbourne dellnitz 93}. We also thank A. Vanderbauwhede for his valuable help when we were in the process of understanding his paper~\cite{vanderbauwhede 95}. Thanks also go to D. Burghelea, M. Field, V. Ginzburg, A. Hern\'andez, K. Hess,  J. E. Marsden, and J. Montaldi for their assistance at various points in the development of this work.

\end{document}